\documentclass[reqno]{amsart}

\usepackage{amssymb}
\usepackage{amsthm}
\usepackage{amsmath}

\usepackage[foot]{amsaddr}
\usepackage[usenames]{color}
\usepackage{hyperref}

\numberwithin{equation}{section}

\theoremstyle{plain}
\newtheorem{proposition}{Proposition}[section]
\newtheorem{theorem}[proposition]{Theorem}
\newtheorem{lemma}[proposition]{Lemma}
\newtheorem{corollary}[proposition]{Corollary}
\newtheorem{result}[proposition]{Result}

\theoremstyle{definition}

\newtheorem{definition}[proposition]{Definition}
\newtheorem{observation}[proposition]{Observation}

\theoremstyle{remark}
\newtheorem{remark}[proposition]{Remark}

\DeclareMathOperator{\Real}{Re}
\DeclareMathOperator{\Imaginary}{Im}

\DeclareMathOperator{\Hol}{Hol}

\DeclareMathOperator{\Euc}{Euc}

\DeclareMathOperator{\Bc}{\mathcal{B}}

\DeclareMathOperator{\Cb}{\mathbb{C}}

\DeclareMathOperator{\Nb}{\mathbb{N}}

\DeclareMathOperator{\Rb}{\mathbb{R}}

\DeclareMathOperator{\Zb}{\mathbb{Z}}

\newcommand{\abs}[1]{\left|#1\right|}
\newcommand{\conj}[1]{\overline{#1}}
\newcommand{\norm}[1]{\left\|#1\right\|}
\newcommand{\wt}[1]{\widetilde{#1}}
\newcommand{\wh}[1]{\widehat{#1}}

\newcommand{\bv}[1]{{#1}^{\bullet}}

\newcommand{\bcdot}{\boldsymbol{\cdot}}
\newcommand{\nrml}{\boldsymbol{\nu}}
\newcommand{\pr}{{\sf p}}
\newcommand{\lemn}{\mathfrak{L}}
\newcommand{\z}{z^{\raisebox{-1pt}{$\scriptstyle {\prime}$}}}
\newcommand{\uni}{{\sf U}}

\begin{document}

\title[Goldilocks domains, visibility, and applications]{Goldilocks domains, a weak notion of visibility, \\ 
and applications}

\author{Gautam Bharali}
\address{Department of Mathematics, Indian Institute of Science, Bangalore 560012, India}
\email{bharali@math.iisc.ernet.in}

\author{Andrew Zimmer}
\address{Department of Mathematics, University of Chicago, Chicago, IL 60637.}
\email{aazimmer@uchicago.edu}

\date{\today}

\begin{abstract} In this paper we introduce a new class of domains in complex Euclidean space, called Goldilocks domains, and study their complex geometry. These domains are defined in terms of a lower bound on how fast the Kobayashi metric grows and an upper bound on how fast the Kobayashi distance grows as one approaches the boundary. Strongly pseudoconvex domains and weakly pseudoconvex domains of finite type always satisfy this Goldilocks condition, but we also present families of Goldilocks domains that have low boundary regularity or have boundary points of  infinite type. We will show that the Kobayashi metric on these domains behaves, in some sense, like a negatively curved Riemannian metric. In particular, it satisfies a visibility condition in the sense of Eberlein and O'Neill. This behavior allows us to prove a variety of results concerning boundary extension of maps and to establish Wolff--Denjoy theorems for a wide collection of domains.
\end{abstract} 

\maketitle

\section{Introduction}\label{sec:intro}

Given a bounded domain $\Omega \subset \Cb^d$, let $K_\Omega: \Omega\times \Omega \rightarrow \Rb_{\geq 0}$ be the Kobayashi distance and let $k_\Omega: \Omega\times \Cb^d \rightarrow \Rb_{\geq 0}$ be the infinitesimal Kobayashi metric (also known as the Kobayashi--Royden metric). One aim of this paper is to present some new techniques\,---\,which we shall use to study the behavior of holomorphic maps into a bounded domain $\Omega$\,---\,that arise from certain intuitions in metric geometry applied to the metric space $(\Omega, K_{\Omega})$. A class of domains that has been studied extensively in the literature is the class of bounded pseudoconvex domains of finite type (which includes the class of bounded strongly pseuodconvex domains). For such a domain $\Omega$, there exist constants $c, \epsilon>0$ such that 
\begin{align*}
k_\Omega(x;v) \geq \frac{c\norm{v}}{\delta_\Omega(x)^\epsilon} \; \; \text{for all $x \in \Omega$ and $v \in \Cb^d$},
\end{align*}
where $\delta_\Omega(x)$ is the distance from $x$ to $\partial\Omega$ (see \cite{C1992}; also see Section~\ref{sec:examples} for explanations). Moreover, in this case, because $\partial\Omega$ has (at least) $C^2$ regularity, it is easy to establish that for each $x_0 \in \Omega$ there exists a constant $C> 0$ (depending on $x_0$) such that 
\begin{align*}
K_\Omega(x_0, x) \leq C + \frac{1}{2}\log \frac{1}{\delta_\Omega(x)}  \; \; \text{for all $x \in \Omega$}.
\end{align*}
We shall show that similar bounds on the growth of $k_\Omega$ and $K_\Omega$ as one approaches the boundary underlie a variety of results on the complex geometry and complex dynamics associated to $\Omega$ (for any bounded domain $\Omega$ admitting such bounds).

To measure the growth rate of the Kobayashi metric as one approaches the boundary we introduce the following function for a bounded domain $\Omega \subset \Cb^d$:
\begin{align*}
M_\Omega(r) := \sup\left\{ \frac{1}{k_{\Omega}(x;v)} : \delta_\Omega(x) \leq r, \norm{v}=1\right\}.
\end{align*}
We are now in a position to define the following class of domains:

\begin{definition}\label{def:good_domains}
A bounded domain $\Omega \subset \Cb^d$ is a \emph{Goldilocks domain} if
\begin{enumerate}
\item for some (hence any) $\epsilon >0$ we have
\begin{align*}
\int_0^\epsilon \frac{1}{r} M_\Omega\left(r\right) dr < \infty,
\end{align*}
\item for each $x_0 \in \Omega$ there exist constants $C, \alpha > 0$ (that depend on $x_0$) such that 
\begin{align*}
K_\Omega(x_0, x) \leq C + \alpha \log \frac{1}{\delta_\Omega(x)}
\end{align*}
for all $x \in \Omega$. 
\end{enumerate}
\end{definition}

\begin{remark}
Given a bounded domain $\Omega$, the function $M_{\Omega}$ is clearly monotone decreasing, hence is Lebesgue measurable. Thus, the integral in (1) of Definition~\ref{def:good_domains} is well defined.
\end{remark}

\begin{remark}\label{rem:just_right}
The slightly unusual phrase ``Goldilocks domain'' is intended to point to the fact that if $\Omega$ is a Goldilocks domain, then $\partial\Omega$ lies in between (and avoids) the extremes of having outward-pointing cusps and having points at which $\partial\Omega$ is flat to infinite order and is, in a precise sense, too flat. A classical argument for planar domains, for instance, implies that the first situation is ruled out by Condition~2 above. Condition~1, it turns out, rules out domains that are not pseudoconvex: i.e., \emph{Goldilocks domain are pseudoconvex}. We discuss all this in more detail in Section~\ref{sec:examples}.
\end{remark}

We have \emph{deliberately} chosen to define Goldilocks domains in rather abstract terms. One objective of this work is to introduce methods that are rooted in metric geometry and are applied to the metric space $(\Omega, K_\Omega)$. The crucial properties that animate these methods are most clearly illustrated when working with domains that satisfy the conditions in Definition~\ref{def:good_domains}.

One deduces from the first paragraph of this section that bounded pseudoconvex domains of finite type are always Goldilocks domains. We shall see, in Section~\ref{sec:examples}, a diverse range of other domains\,---\,described in more geometrically explicit terms\,---\,that are Goldilocks domains. Consequently, we are able to establish, among several other results, extensions of the following widely-studied phenomena:
\begin{itemize}
 \item Wolff--Denjoy theorems in higher dimensions,
 \item continuous boundary-extension of proper holomorphic maps
\end{itemize}
to a new range of domains. To give a sense of the uses that the methods hinted at are put to: a pseudoconvex domain $\Omega\Subset \Cb^d$, $d\geq 2$, of finite type is, in general, non-convex and there may be points $\xi \in \partial\Omega$ around which $\Omega$ is not locally convexifiable. Such an $\Omega$ has \emph{very little} resemblance to a convex domain, yet a form of the Wolff--Denjoy theorem holds on $\Omega$\,---\,see Corollary~\ref{cor:WD_finite-type}.

We now introduce the main theorems of this paper.

\subsection{Negative curvature and the Kobayashi metric} It is well known that the unit ball $\Bc \subset \Cb^d$ endowed with the Kobayashi metric is isometric to complex hyperbolic space, which is a canonical example of a negatively curved Riemannian manifold. Based on this example, it is natural to conjecture that domains that are similar to the unit ball also have a negatively curved Kobayashi metric. One problem with this conjecture is that the infinitesimal Kobayashi metric on a general domain is not Riemannian and has low regularity, making a notion of infinitesimal curvature difficult to define. One remedy to this problem is to consider a coarse notion of negative curvature introduced by Gromov~\cite{G1987} which is now called Gromov hyperbolicity. 

Along these lines Balogh and Bonk~\cite{BB2000} proved that the Kobayashi distance on a bounded strongly pseudoconvex domain is Gromov hyperbolic. The Kobayashi distance is also Gromov hyperbolic for bounded convex domains whose boundary has finite type in the sense of D'Angelo~\cite{Z2014}. 

We will show, in Section~\ref{sec:examples}, some examples of Goldilocks domains where the Kobayashi distance is not Gromov hyperbolic. However, a key part of this paper is to show that the Kobayashi metric on a Goldilocks domain does have some negative-curvature-like behavior. In particular we are motivated by a definition of Eberlein and O'Neill~\cite{EO1973} who call a non-positively curved simply connected Riemannian manifold $X$ \emph{a visibility space}  if for every two points $\xi, \eta$ in the ideal boundary $\partial X$ and neighborhoods $V_\xi, V_\eta$ of $\xi,\eta$ in $X \cup \partial X$ so that $\overline{V_\xi} \cap \overline{V_\eta} = \emptyset$ there exists a compact set $K \subset X$ with the following property: if $\sigma: [0,T] \rightarrow X$ is a geodesic with $\sigma(0) \in V_\xi$ and $\sigma(T) \in V_\eta$ then $\sigma \cap K \neq \emptyset$ (see also ~\cite[page 54]{BGS1985} or~\cite[page 294]{BH1999}). Informally, this states that geodesics between two distinct points at infinity bend into the space. 

It is well known that a complete negatively curved simply connected Riemannian manifold $X$ is always a visibility space and, more generally, a proper geodesic Gromov hyperbolic metric space always satisfies a visibility type condition (see for instance~\cite[page 428]{BH1999}). 

In the context of a Goldilocks domain $\Omega$, we do not know that the metric space $(\Omega, K_\Omega)$ is Cauchy complete and in particular we do not know whether or not every two points can be joined by a geodesic. This leads us to consider a more general class of curves which we call \emph{almost-geodesics} (defined in Section~\ref{sec:curves}). We will then prove: 

\begin{theorem}\label{thm:visible} (see Section~\ref{sec:visible})
Suppose $\Omega \subset \Cb^d$ is a Goldilocks domain and $\lambda \geq 1$, $\kappa \geq 0$. If $\xi,\eta \in \partial\Omega$ and $V_\xi, V_\eta$ are neighborhoods of $\xi,\eta$ in $\overline{\Omega}$ so that $\overline{V_\xi} \cap \overline{V_\eta} = \emptyset$ then there exists a compact set $K \subset \Omega$ with the following property: if $\sigma: [0,T] \rightarrow \Omega$ is an $(\lambda, \kappa)$-almost-geodesic with $\sigma(0) \in V_\xi$ and $\sigma(T) \in V_\eta$  then $\sigma \cap K \neq \emptyset$. 
\end{theorem}

This theorem makes intuitive sense: on a Goldilocks domain the Kobayashi metric grows rapidly as one approaches the boundary and so length minimizing curves wish to spend as little time as possible near the boundary. This leads to the phenomenon of such curves bending into the domain and intersecting some fixed compact set. A key point of this paper is giving a precise condition on the rate of blow-up, namely Definition~\ref{def:good_domains}, which leads to this behavior. 

There are several visibility type results in the literature. Chang, Hu, and Lee studied the limits of complex geodesics in strongly convex domains and proved a visibility type result for complex geodesics, see~\cite[Section 2]{CHL1988}. Mercer~\cite{M1993} extended these results to \emph{$m$-convex domains}, that is bounded convex domains $\Omega$ for which there exist constants $C > 0$ and $m > 0$ such that
\begin{equation}
\label{est:m_convex}
\inf\left\{ \norm{x-\xi} : \xi \in \partial \Omega \cap (x + \Cb \bcdot v)\right\} \leq C\delta_{\Omega}(x)^{1/m} \; \;
\forall x\in \Omega \text{ and } \forall v\in \Cb^d\setminus\{0\}.
\end{equation}
Notice that every strongly convex set is $2$-convex. Finally, Karlsson~\cite[Lemma 36]{K2005b} proved a visibility result for geodesics in bounded domains $\Omega$ that satisfy estimate~\eqref{est:m_convex}, have Cauchy-complete Kobayashi metric, and have $C^{1,\alpha}$ boundary.

\subsection{Continuous extensions of proper holomorphic maps}
The earliest result on the continuous extension up to the boundary of a proper holomorphic map between a pair of domains $D$ and $\Omega$ in $\Cb^d$, $d\geq 2$, with no other assumptions on the map, was established by Pin{\v{c}}uk \cite{P1974}. In that work, the domains $D$ and $\Omega$ are assumed to be strongly pseudoconvex with $C^2$ boundaries. Owing to strong pseudoconvexity, it is shown that the continuous extension of the map from $D$ to $\Omega$ satisfies a H{\"o}lder condition on $\overline{D}$ with H{\"o}lder exponent $1/2$. Soon thereafter, the focus of the question of boundary regularity of a proper holomorphic map between a pair of domains of the same dimension shifted largely to domains with $C^\infty$ boundaries, and to obtaining \emph{smooth} extension to the boundary, of the given proper map, largely due to various Bergman-kernel methods introduced by Fefferman \cite{F1974}\,---\,who considered biholomorphisms\,---\,and by Bell--Ligocka \cite{BL1980} and Bell \cite{B1981}. The literature on the smooth extension of proper holomorphic maps is truly enormous, and we refer the interested reader to the survey \cite{F1993} by Forstneri{\v{c}}.

The latter methods are not helpful if the boundary of either of $D$ or $\Omega$ has low regularity, or if either domain is assumed to be merely pseudoconvex (i.e., without any finite-type condition on the boundary). In this context, the methods of Diederich--Forn{\ae}ss \cite{DF1979} are helpful. The idea of using the Kobayashi metric was first introduced in \cite[Theorem~1]{DF1979}, and, in that theorem, the requirement of strong pseudoconvexity of $D$ in Pin{\v{c}}uk's theorem is dropped. In this paper, we generalize \cite[Theorem~1]{DF1979} by allowing the target domain to have non-smooth boundary. We point to Section~\ref{sec:examples} for a sense of how irregular $\partial\Omega$, for $\Omega$ as in Theorem~\ref{thm:proper} below, can be. One frequently encounters proper holomorphic maps of smoothly-bounded domains whose images have non-smooth boundary; consider the various examples of maps between Reinhardt domains. In the latter setting and in low dimensions, continuous extension up to the boundary follows from an explicit description of the proper map in question\,---\,see \cite{IK2006}, for instance. Even in $\Cb^2$ though, establishing such descriptions as in \cite{IK2006} is a highly technical effort. In contrast, the question of continuous extension\,---\,and for a variety of boundary geometries for the target space\,---\,is settled by the following theorem.    
 
 \begin{theorem}\label{thm:proper} (see Section~\ref{sec:proper})
 Let $D$ and $\Omega$ be bounded domains in $\Cb^d$. Suppose $D$ is pseudoconvex with $C^2$-smooth boundary, and $\Omega$ is a Goldilocks domain satisfying an interior-cone condition. Any proper holomorphic map $F: D\rightarrow \Omega$ extends to a continuous map on $\conj{D}$.
\end{theorem}

We refer the reader to Section~\ref{sec:examples} for a definition of the interior-cone condition. This cone condition on $\Omega$ above  allows us to adapt certain ideas in \cite{DF1979}. Here is a sketch of the proof: using a type of Hopf lemma, which is available due to the cone condition on $\Omega$,  we first show that $\delta_\Omega(F(z)) \leq c \delta_D(z)^{\eta}$ for some $c, \eta > 0$. Now suppose $\xi \in \partial D$ and $\nrml(\xi)$ is the inward-pointing normal ray, then the rapid growth of the Kobayashi metric is used to show that the curve $F(\xi+t\nrml(\xi))$ does not oscillate very much as $t \searrow 0$ and, in particular, one obtains a continuous extension to $\partial D$. As in the Theorem~\ref{thm:visible} above, the key point is to have the precise rate of blow-up necessary to obtain such behavior.

\subsection{Continuous extensions of quasi-isometric embeddings}\label{ssec:cont_extn} 

We can also prove continuous extensions of certain non-holomorphic maps between domains of different dimensions. A map $F: (X, d_X) \rightarrow (Y, d_Y)$ between two metric spaces is called a \emph{$(\lambda, \kappa)$-quasi-isometric embedding} if there exist constants $\lambda \geq 1$ and $\kappa \geq 0$ so that
\begin{align*}
 \frac{1}{\lambda} d_Y(F(x_1), F(x_2))-\kappa  \leq d_X(x_1, x_2) \leq \lambda d_Y(F(x_1), F(x_2))+\kappa
\end{align*}
for all  $x_1,x_2\in X$.

There are two motivations for investigating continuous extensions of quasi-isometric embeddings. Our main motivation stems from Lempert's theorem \cite{L1981}\,---\,and its generalization to convex domains with non-smooth boundaries by Royden and Wong \cite{RW1983}\,---\,which establish that there exist complex geodesics between any pair of points of a convex domain in $\Cb^d$, $d\geq 2$. A \emph{complex geodesic} of a domain $\Omega$ is a holomorphic map from $\Delta$ (the open unit disk in $\Cb$) to $\Omega$ that is an isometric embedding of $(\Delta, K_{\Delta})$ into $(\Omega, K_{\Omega})$. It is natural to ask whether a complex geodesic extends continuously up to $\partial\Delta$. This question has been examined\,---\,beginning with Lempert's result for strongly convex domains with $C^3$-smooth boundary \cite{L1981}\,---\,from various perspectives \cite{M1993, B2016, Z2015}, but:
\begin{itemize}
 \item a comprehensive answer to this question is still forthcoming;
 \item little is known in general, at present, for domains that are \emph{non-convex} and admit complex geodesics. 
\end{itemize}

We ought to mention that all complex geodesics of the symmetric twofold product of $\Delta$ (also known as the \emph{symmetrized bidisc})\,---\,which is not biholomorphic to any convex domain; see \cite{C2004}\,---\,extend continuously up to $\partial\Delta$. This follows from the work of Agler--Young \cite{AY2004} and Pflug--Zwonek \cite{PZ2005}. Those results rely heavily on the specific properties of the symmetrized bidisc. 

However, there is a general approach to answering this question. When $(X, d_X)$ is a proper geodesic Gromov hyperbolic metric space, $X$ has a natural boundary $X(\infty)$ ``at infinity'', and the set $X\cup X(\infty)$ has a topology that makes it a compactification of $X$ (see Section~\ref{sec:gromov_prod} below for more details). One of the fundamental properties of this compactification is the extension of quasi-isometries:

\begin{result} \label{res:gromov_qi_ext}(see for instance~\cite[Chapter III.H, Theorem 3.9]{BH1999})
Suppose $(X,d_X)$ and $(Y,d_Y)$ are two proper geodesic Gromov hyperbolic metric spaces. Then any continuous quasi-isometric embedding $F:(X, d_X) \rightarrow (Y,d_Y)$ extends to a continuous map $\wt{F} :X \cup X(\infty) \rightarrow Y \cup Y(\infty)$. 
\end{result}

It is very easy to see that $(\Delta, K_{\Delta})$ is Gromov hyperbolic, with $\Delta(\infty) = \partial\Delta$. Thus, if one could show that $(\Omega, K_{\Omega})$ satisfies all the conditions in Result~\ref{res:gromov_qi_ext} and that $\Omega(\infty) = \partial\Omega$\,---\,where $\Omega$ is a domain that admits complex geodesics\,---\,then one would have an answer to the above question.

However, by the main theorem of \cite{Z2014}, if $\Omega\subset \Cb^d$, $d\geq 2$, is a smoothly bounded convex domain having infinite-type points $\xi\in \partial\Omega$ (i.e., $T_\xi(\partial\Omega)$ has infinite order of contact with $\partial\Omega$ along a \emph{complex} direction in $T_\xi(\partial\Omega)$) then $(\Omega, K_{\Omega})$ is not Gromov hyperbolic. Thus, approaches other than Result~\ref{res:gromov_qi_ext} are of interest. Independently of all this, it would be interesting in itself to prove an analogue of Result~\ref{res:gromov_qi_ext} in which, working in the category of domains in $\Cb^d$, the Gromov-hyperbolicity assumption on either of $(X, d_X)$ or $(Y, d_Y)$ (or both) is supplanted by a \emph{strictly weaker} assumption by taking advantage of our knowledge of the Kobayashi metric. The latter is further motivation for the following analogue of Result~\ref{res:gromov_qi_ext}:
 
\begin{theorem}(see Theorem~\ref{thm:quasi_isometry_ext} below)\label{thm:qi_ext}
Let $D$ be a bounded domain in $\Cb^k$ and suppose $(D, K_D)$ is a proper geodesic Gromov hyperbolic metric space. Let $\Omega\subset \Cb^d$ be a Goldilocks domain. If $F :(D, K_D) \rightarrow (\Omega, K_\Omega)$ is a continuous quasi-isometric embedding, then there exists a continuous extension $\wt{F} : D \cup D(\infty) \rightarrow \overline{\Omega}$.
\end{theorem}
 
\begin{remark} \ 
\begin{enumerate}
\item Our proof of Theorem~\ref{thm:qi_ext} will follow from that of Theorem~\ref{thm:quasi_isometry_ext} below.  Theorem~\ref{thm:quasi_isometry_ext} is \emph{much more general} and the techniques used in its proof apply to a wide range of metric spaces (even with $(D, K_D)$ replaced by more general metric spaces) and compactifications. However, we shall focus on domains in this paper.
\item Theorem~\ref{thm:qi_ext} and Theorem~\ref{thm:quasi_isometry_ext} both represent applications of visibility. A key step of the proof is Proposition~\ref{prop:qg_visibility} below, where we establish a visibility result for quasi-geodesics.
\item If $\Omega$ is strongly pseudoconvex or convex with finite-type boundary, then $(\Omega,K_\Omega)$ is a proper geodesic Gromov hyperbolic metric space: see \cite{BB2000} and \cite{Z2014}, respectively. Hence in these cases,  Theorem~\ref{thm:qi_ext} follows directly from Result~\ref{res:gromov_qi_ext}. However, proving that the Kobayashi metric is Gromov hyperbolic in either case is very involved, and our approach of using a visibility condition is much more direct.  
\end{enumerate}
 \end{remark}

Going back to our initial motivation for Theorem~\ref{thm:qi_ext}: it follows from this theorem that if $\varphi: \Delta \rightarrow \Omega$ is a complex geodesic  into a Goldilocks domain, then $\varphi$ extends to a continuous map $\wt{\varphi} : \overline{\Delta} \rightarrow \overline{\Omega}$. This, in fact, extends known results to the case when $\Omega$ is not necessarily convex. We refer the reader to subsection~\ref{ssec:implications} below.

\subsection{Wolff--Denjoy Theorems}\label{ssec_WD} There has been considerable interest in understanding the behavior of iterates of a holomorphic map $f:\Omega \rightarrow \Omega$ on a bounded domain $\Omega$. Since $\Omega$ is bounded, for any subsequence $n_i \rightarrow \infty$ one can always find a subsequence $n_{i_j} \rightarrow \infty$ so that $f^{n_{i_j}}$ converges locally uniformly to a holomorphic map $F:\Omega \rightarrow \overline{\Omega}$. The general goal is to show that the behavior of each convergent subsequence is identical. This is demonstrated in the classical Wolff--Denjoy theorem:

\begin{result}[\cite{D1926, W1926}]\label{res:WD}
Suppose $f:\Delta \rightarrow \Delta$ is a holomorphic map then either:
\begin{enumerate}
\item $f$ has a fixed point in $\Delta$; or
\item there exists a point $\xi \in \partial \Delta$ so that
\begin{equation*}
\lim_{n \rightarrow \infty} f^n(x) = \xi
\end{equation*}
 for any $x \in \Delta$, this convergence being uniform on compact subsets of $\Delta$.
\end{enumerate}
\end{result}

The above result was extended to the unit (Euclidean) ball in $\Cb^d$, for all $d$, by Herv{\'e} \cite{H1963}.   
It was further generalized by Abate\,---\,see \cite{A1988} or \cite[Chapter~4]{A1989}\,---\,to strongly convex domains. The above theorem was later generalized to contractible strongly pseudoconvex domains by Hua~\cite{H1984} and to a variety of different types of convex domains (see for instance~\cite{A2014} and the references therein). Wolff--Denjoy theorems are also known to hold on certain metric spaces where a boundary at infinity replaces the topological boundary, see for instance~\cite{K2001} or~\cite{B1997}.

Using the visibility result, we will prove two Wolff--Denjoy theorems for Goldilocks domains. The first theorem concerns holomorphic maps on taut Goldilocks domains while the second theorem considers maps that are 1-Lipschitz with respect to the Kobayashi distance and Goldilocks domains $\Omega$ for which $(\Omega, K_{\Omega})$ Cauchy complete. Since every holomorphic map is 1-Lipschitz with respect to the Kobayashi distance, our second theorem considers a more general class of maps. On the other hand, because whenever $(\Omega, K_{\Omega})$ is Cauchy complete the domain $\Omega$ is taut, our first theorem considers more a general class of domains.

It is not hard to see that the dichotomy presented by Result~\ref{res:WD} fails in general if the domain in question is not contractible. The following theorems present the dichotomy relevant to more general circumstances. Here are the precise statements:

\begin{theorem}\label{thm:WD} (see Section~\ref{sec:WD} below)
Suppose $\Omega \subset \Cb^d$ is a taut Goldilocks domain. If $f:\Omega \rightarrow \Omega$ is a holomorphic map then either:
\begin{enumerate}
\item for any $x \in \Omega$ the orbit $\{ f^n(x): n \in \Nb\}$ is relatively compact in $\Omega$; or
\item there exists $\xi \in \partial \Omega$ so that
\begin{equation*}
\lim_{n \rightarrow \infty} f^n(x) = \xi
\end{equation*}
 for any $x \in \Omega$, this convergence being uniform on compact subsets of $\Omega$.
\end{enumerate}
\end{theorem}

\begin{theorem}\label{thm:m_WD} (see Section~\ref{sec:WD} below)
Suppose $\Omega \subset \Cb^d$ is a  Goldilocks domain such that $(\Omega, K_\Omega)$ is Cauchy complete. If $f:\Omega \rightarrow \Omega$ is 1-Lipschitz with respect to the Kobayashi distance then either:
\begin{enumerate}
\item for any $x \in \Omega$ the orbit $\{ f^n(x): n \in \Nb\}$ is relatively compact in $\Omega$; or
\item there exists $\xi \in \partial \Omega$ so that
\begin{equation*}
\lim_{n \rightarrow \infty} f^n(x) = \xi
\end{equation*}
 for any $x \in \Omega$, this convergence being uniform on compact subsets of $\Omega$.
\end{enumerate}
\end{theorem}

The assumption of Cauchy completeness of $(\Omega, K_{\Omega})$ provides tools, namely work of Ca{\l}ka \cite{C1984b}, that do not have analogues in the taut setting. In particular, the proof of Theorem~\ref{thm:WD} is much more intricate than the proof of Theorem~\ref{thm:m_WD}. However, tautness is a rather mild condition: for instance a bounded pseudoconvex domain with $C^1$ boundary is known to be taut~\cite{KR1981} (whereas it is unknown whether $(\Omega, K_{\Omega})$ is Cauchy complete if $\Omega$ is a weakly pseudoconvex domain of finite type in $\Cb^d$, $d > 2$). This allows one to state various types of corollaries of Theorem~\ref{thm:WD}\,---\,for instance, see Corollary~\ref{cor:WD_finite-type} below.

\subsection{Basic notations} We end the introduction by fixing some very basic  notations. 
\begin{enumerate}
\item For $z \in\Cb^d$, $\norm{z}$ will denote the standard Euclidean norm and, for $z_1, z_2 \in \Cb^d$, $d_{\Euc}(z_1, z_2) = \norm{z_1 - z_2}$ will denote the standard Euclidean distance.
\item $\Delta \subset \Cb$ will denote the open unit disk, and $\rho_\Delta$ will denote the Poincar{\'e} metric on $\Delta$.
\item For a point $z\in \Cb^d$ and $r > 0$, $B_r(z)$ will denote the open Euclidean ball with center $z$ and radius $r$. 
\end{enumerate}

\subsection*{Acknowledgments} 
Gautam Bharali is supported in part by a Swarnajayanti Fellowship (Grant No.~DST/SJF/MSA-02/2013-14) and by a UGC Centre for Advanced Study grant. Andrew Zimmer is partially supported by the National Science Foundation under Grant No.~NSF 1400919.
 
\section{Examples and corollaries}\label{sec:examples}

In this section we shall present certain broad classes of bounded domains\,---\,described in terms of rather explicit boundary properties\,---\,under which either Condition~1 or Condition~2 in the definition of a Goldilocks domain (i.e., Definition~\ref{def:good_domains}) is satisfied. Consequently, we shall see that Definition~\ref{def:good_domains} admits a truly wide range of bounded domains.

\subsection{Domains that satisfy Condition~2}\label{ssec:cond_2}
Lemma~\ref{lem:int_cone} below establishes that a simple property, which arises in several areas in analysis, ensures that any domain with this property satisfies Condition~2. We require a couple of definitions.

\begin{definition}
An \emph{open right circular cone with aperture $\theta$} is an open subset of $\Cb^d$ of the form
\begin{align*}
 \{z\in \Cb^d : \Real[\,\langle z, v\rangle\,] > \cos(\theta/2)\norm{z}\}
 =: \Gamma(v, \theta),
\end{align*}
where $v$ is some unit vector in $\Cb^d$, $\theta\in (0, \pi)$, and $\langle\bcdot\,,\,\bcdot\rangle$ is the standard
Hermitian inner product on $\Cb^d$. For any point $p\in \Cb^d$, the \emph{axis} of the (translated) cone $p+\Gamma(v, \theta)$ is the ray $\{p + tv: t > 0\}$.
\end{definition}

\begin{definition}\label{def:cone_cond}
Let $\Omega$ be a bounded domain in $\Cb^d$. We say that $\Omega$ satisfies an \emph{interior-cone condition
with aperture $\theta$} if there exist constants  $r^0 > 0$, $\theta\in (0, \pi)$, and a compact subset $K\subset \Omega$ such that for each $x\in \Omega\setminus K$, there exist a point $\xi_x\in \partial\Omega$ and a unit vector $v_x$ such that
\begin{itemize}
 \item $x$ lies on the axis of the cone $\xi_x+\Gamma(v_x, \theta)$, and
 \item $(\xi_x+\Gamma(v_x, \theta))\cap B_{r^0}(\xi_x) \subset \Omega$.
\end{itemize}
We say that $\Omega$ \emph{satisfies an interior-cone condition} if there exists a $\theta \in (0, \pi)$ so that $\Omega$ satisfies an interior-cone condition with aperture $\theta$.
\end{definition}

The proof of the following statement involves a mild adaptation of a technique used in \cite[Proposition~2.5]{FR1987} and in \cite[Proposition~2.3]{M1993}.

\begin{lemma}\label{lem:int_cone}
Let $\Omega$ be a bounded domain in $\Cb^d$ that satisfies an interior-cone condition with aperture
$\theta$. Then $\Omega$ satisfies Condition~2 in the definition of a Goldilocks domain.
\end{lemma}

\begin{proof}
For any $\beta > 1$, define the holomorphic map $\psi_\beta: \Delta\rightarrow \Cb$ by
\begin{align*}
 \psi_\beta(\zeta) := (1+\zeta)^{1/\beta}.
\end{align*}
Given a unit vector $v\in \Cb^d$ and a number $r > 0$, define the holomorphic map
$\Psi(\bcdot\,;\,\beta, v, r): \Delta\rightarrow \Cb^d$ by
\begin{align*}
 \Psi(\zeta; \beta, v, r) := r\psi_\beta(\zeta)v,
\end{align*}
and denote the image of  $\Psi(\bcdot\,;\,\beta, v, r)$ by $\lemn(\beta, v, r)$.

It is an elementary calculation that there exist constants $R > 0$ and $\alpha > 1$ such that
\begin{align*}
 R\psi_\alpha(\Delta) \subset \{\zeta\in \Cb: \Real(\zeta) > \cos(\theta/2)|\zeta|\}\cap \{\zeta\in \Cb : |\zeta| < r^0\},
\end{align*}
where $\theta$ and $r^0$ are as given by Definition~\ref{def:cone_cond}. It follows from this, owing to our condition on $\Omega$, that:
\begin{itemize}
 \item[$(\bullet)$] There exists a compact subset $K^\prime$ such that $K\subset K^\prime\subset \Omega$ and such that for each $x\in \Omega\setminus K^\prime$, there exist a point $\xi_x\in \partial\Omega$ and a unit vector $v_x$ so that
 \begin{itemize}
  \item[$(i)$] $\xi_x + \lemn(\alpha, v_x, R)\subset \Omega$;
  \item[$(ii)$] $x$ lies on the line segment joining $\xi_x$ to $\xi_x+\Psi(0; \alpha, v, R) =: q_x$; and
  \item[$(iii)$] $q_x\in K^\prime$.
 \end{itemize}
\end{itemize}
Then, for $x\in \Omega\setminus K^\prime$, there exists a unique number $t(x) > 0$ such that $\xi_x + t(x)v_x = x$. Clearly $\delta_\Omega(x)\leq
t(x)$. Also, $\Psi(\bcdot\,;\,\alpha, v, R)$ maps the point
\begin{align*}
 \big((t(x)/R)^\alpha - 1\big) \in (-1, 0)
\end{align*}
to the point $x$.

Fix $x_0\in \Omega$. It suffices to establish the inequality that defines Condition~2 for $x\in K^\prime$. Set
$C_1 := \sup\{K_\Omega(z, x_0)  : z\in K^\prime\}$. Then, by $(\bullet)$, if $x\in \Omega\setminus K^\prime$,
then
\begin{align*}
 K_\Omega(x_0, x) \leq K_\Omega(x_0, q_x) + K_\Omega(q_x, x) &\leq C_1 + \rho_\Delta(0, (t(x)/R)^\alpha - 1) \\
 &= C_1 + \rho_\Delta(0, 1- (t(x)/R)^\alpha ) \\
 &\leq C_1 + \frac{1}{2}\log\left(\frac{2}{(t(x)/R)^\alpha}\right) \\
 &\leq \big(C_1 + (1/2)\log(2R^\alpha)\big) + (\alpha/2)\log\frac{1}{\delta_\Omega(x)}.
\end{align*}
Hence, $\Omega$ satisfies Condition~2.
\end{proof}

This gives us the following:

\begin{corollary}
Let $\Omega_1$ and $\Omega_2$ be two convex Goldilocks domains in $\Cb^d$ having non-empty intersection. Then
$\Omega_1\cap \Omega_2$ is also a Goldilocks domain.
\end{corollary}
\begin{proof}
Write $D = \Omega_1\cap \Omega_2$. Since $D$ is a convex domain, it satisfies an interior-cone condition with aperture $\theta$ for some $\theta\in (0, \pi)$. Thus, by Lemma~\ref{lem:int_cone}, $D$ satisfies Condition~2.

Since $D\subset \Omega_j$, $j = 1, 2$, we have
\begin{equation}\label{eq:k_D_bigger}
 k_D(x; v) \geq k_{\Omega_j}(x; v) \; \; \forall x\in D, \ \forall v: \|v\| = 1, \text{ and } j = 1, 2.
\end{equation}
Fix an $r > 0$. Then\vspace{-2mm}
\begin{multline*}
 \{x\in  D: \delta_D(x)\leq r\} \\
 \subseteq
 \{x\in D : \delta_{\Omega_1}(x)\leq r\}\cup \{x\in D : \delta_{\Omega_2}(x)\leq r\}\,\equiv\,\mathcal{S}(1, r)
 \cup \mathcal{S}(2, r).
\end{multline*}
Thus, by \eqref{eq:k_D_bigger}, we can estimate:
\begin{align*}
 M_D(r) &\leq \sup_{(\mathcal{S}(1, r)\cup \mathcal{S}(2, r))\times \{\norm{v} = 1\}}\frac{1}{k_D(x; v)} \\
 &= \max\Big[\sup_{\mathcal{S}(1, r)\times \{\norm{v} = 1\}}
       \frac{1}{k_D(x; v)}\,,\;\sup_{\mathcal{S}(2, r)\times \{\norm{v} = 1\}}\frac{1}{k_D(x; v)}\Big] \\
 &\leq \max\big(M_{\Omega_1}(r), M_{\Omega_2}(r)\big).
\end{align*}
Now, $M_D$, being monotone increasing, is Lebesgue measurable. Since $M_{\Omega_1}$ and $M_{\Omega_2}$ satisfy the inequality that defines Condition~1, the above estimate ensures that $M_D$ does so too. Thus, $D$ satisfies Condition~1. Hence, $D$ is a Goldilocks domain.
\end{proof}

\subsection{Domains that satisfy Condition~1}\label{ssec:cond_1}
 In looking for domains that satisfy Condition~1, we shall examine two classes of domains with very different degrees of boundary smoothness. Let us first examine a class of domains with $C^\infty$-smooth boundaries. In this connection, we need the following result.

\begin{result}[Cho, \cite{C1992}]
Let $\Omega$ be a bounded domain in $\Cb^d$, let $\partial\Omega\cap U$ be smooth and pseudoconvex, where $U$ is a neighborhood of a point $\xi_0\in \partial\Omega$, and let $\partial\Omega$ be of finite 1-type in the sense of D'Angelo at $\xi_0$. Then there exist a neighborhood $V\subset U$ of $\xi_0$ and constants $c, \epsilon > 0$ such that for every $z\in \Omega\cap V$ and for every $v\in \Cb^d$,
\begin{align*}
 k_{\Omega}(z; v) \geq \frac{c\norm{v}}{\delta_{\Omega}(z)^{\epsilon}}.
\end{align*}
\end{result}

The following is now straightforward.

\begin{lemma}
Let $\Omega$ be a bounded pseudoconvex domain of finite type. Then $\Omega$ satisfies Condition~1 in the definition of a Goldilocks domain.
\end{lemma}

\begin{proof}
By the above result, and owing to our hypothesis, we can find finitely many connected open sets $V_1,\dots, V_N$ that cover $\partial\Omega$ and constants $\epsilon_1,\dots, \epsilon_N$ such that
\begin{align*}
  k_{\Omega}(z; v) \geq c\delta_{\Omega}(z)^{-\epsilon_j}
\end{align*}
for every $z\in \Omega\cap V_j$ and for every unit vector $v$, where $c > 0$ is a suitable constant. Write
$s := \min(\epsilon_1,\dots, \epsilon_N)$. Then, for $r > 0$ so small that $r < 1$ and 
\begin{align*}
 \{z\in \Omega : \delta_{\Omega}(z)\leq r\} \subset V_1\cup\dots\cup V_N,
\end{align*}
we have $M_{\Omega}(r)\leq (1/c)r^{s}$, $s > 0$, whence Condition~1 is satisfied.
\end{proof}

The second family of domains that we shall consider will be bounded convex domains. As has been emphasized in Section~\ref{sec:intro}, we would like to consider domains $\Omega$ such that, at any smooth point $\xi\in \partial\Omega$, the boundary is allowed to osculate $H_\xi(\partial\Omega) := T_\xi(\partial\Omega)\cap iT_\xi(\partial\Omega)$ to infinite order and, yet, are not necessarily smoothly bounded. One needs a device to quantify how flat $\partial\Omega$ can get at smooth points. This is accomplished by the notion of the {\em support of $\Omega$ from the outside}, which was introduced in \cite{B2016}. The following definition has been adapted from \cite{B2016}\,---\,which focuses on domains with $C^1$-smooth boundary\,---\,to admit convex domains with non-smooth boundaries as well. (Augmenting our notation somewhat, we shall write $B^{k}_r(z)$ to denote the open Euclidean
ball in $\Cb^k$ with center $z$ and radius $r$.)

\begin{definition}\label{def:supp}
Let $\Omega$ be a bounded convex domain in $\Cb^d, \ d\geq 2$. Let $F: B^{d-1}_r(0)\rightarrow \Rb$ be a $C^1$-smooth convex function with $F(0)=0$ and $DF(0)=0$. We say that {\em $F$ supports $\Omega$ from the outside} if there exists a constant $R\in (0,r)$ such that, for each point $\xi\in \partial\Omega$, there exists a unitary transformation $\uni_\xi$ so that
\begin{itemize}
 \item the set $\big(\xi+\uni_\xi^{-1}(\{v\in \Cb^d: v_d=0\})\big)$ is a supporting complex hyperplane of $\Omega$ at $\xi$, and
 \item the line $\big(\xi+\uni_\xi^{-1}({\rm span}_{\Rb}\{(0,\dots,0,i)\})\big)$ intersects $\Omega$,
\end{itemize}
and such that, denoting the $\Cb$-affine map $v\!\longmapsto\!\uni_\xi(v-\xi)$ as $\uni^\xi$, we have\vspace{-1mm}
\begin{equation*}
 \uni^\xi(\overline\Omega)\cap \big(B^{d-1}_R(0)\times\Delta\big)\,\subset\,\{z=(\z, z_d)\in B^{d-1}_R(0)\times\Delta
 : \Imaginary(z_d) \geq F(\z)\}.
\end{equation*}
\end{definition}

This notion allows us to describe another family of domains that satisfy Condition~1. However, to do so, we will need the following result.

\begin{result}[Graham, \cite{G1990, G1991}]\label{res:kob_bounds}
Let $\Omega$ be a bounded convex domain in $\Cb^d$. For each $z\in \Omega$
and $v\in \Cb^d\setminus\{0\}$, define
\begin{align*}
 r_{\Omega}(z; v) := \sup\Big\{r > 0: \Big(z+ (r\Delta)\frac{v}{\norm{v}}\,\Big)\subset \Omega\Big\},
\end{align*} 
Then:
\begin{equation}\label{eq:kob_bounds}
 \frac{\norm{v}}{2r_{\Omega}(z; v)}\,\leq\,k_{\Omega}(z; v)\,\leq\,\frac{\norm{v}}{r_{\Omega}(z; v)} \quad
 \forall z\in \Omega  \text{ and }  \forall v\in \Cb^d\setminus\{0\}.
\end{equation}
\end{result}

\noindent{The lower bound on $k_\Omega(z;v)$ is the non-trivial part of the result and a proof can also be found in \cite[Theorem 2.2]{F1991}.}

\begin{lemma}\label{lem:convex_cond_1}
Let $\Omega$ be a bounded convex domain in $\Cb^d$, $d\geq 2$. Let $\varPsi : [0, r)\rightarrow \Rb$ be a convex, strictly increasing $C^1$ function such that
\begin{align*}
 \int\nolimits_0^{\epsilon}t^{-1}\varPsi^{-1}(t)\,dt < \infty
\end{align*}
(where $\epsilon > 0$ is small enough for $\varPsi^{-1}$ to be defined). Assume that $\Omega$ is supported from the outside by $F(\z) := \varPsi(\,\norm{\z}\,)$ (write $z = (\z, z_d)$ for each $z\in \Cb^d$). Then $\Omega$ satisfies Condition~1 in the definition of a Goldilocks domain.
\end{lemma}

\begin{proof}
Let $R$ be as given by Definition~\ref{def:supp} with $F = \varPsi(\,\norm{\bcdot}\,)$. Let $C:= \sup_{t\in [0, R)}\varPsi(t)$ and define $t_0$ as follows (it is easy to argue that the set on the right is finite):
\begin{equation}\label{eq:cap}
 t_0 := \min\left[\{C/2\}\cup \{t\in (0,C) : t = \varPsi^{-1}(t)\}\right].
\end{equation}
Let us define
\begin{align*}
 \boldsymbol{{\sf M}} := \{(\z, z_d)\in B^{d-1}_r(0)\times\Cb : \Imaginary(z_n) = \varPsi(\,\norm{\z}\,)\}.
\end{align*}
Let $z\in \Omega$ and let $\xi(z)\in \partial\Omega$ be such that $\delta_{\Omega}(z) = d_{\Euc}(z, \xi(z))$. Clearly
\begin{align*}
 \xi(z)\in \partial B_{\delta_{\Omega}(z)}(z), \; \; \; \text{and} \; \; \;
 \partial\Omega\cap B_{\delta_{\Omega}(z)}(z) = \emptyset,
\end{align*}
whence $B_{\delta_{\Omega}(z)}(z)\subset \Omega$. Thus, for any $(d-1)$-dimensional complex subspace $E$ such that $E\neq H_{\xi(z)}(\partial B_{\delta_{\Omega}(z)}(z))$, the $\Cb$-affine subspace $(\xi(z)+E)$ intersects $B_{\delta_{\Omega}(z)}(z)$, and hence intersects $\Omega$. Therefore, at any point $\xi\in \partial\Omega$ that is of the form $\xi(z)$ for some $z\in \Omega$, there is a unique supporting complex hyperplane of $\Omega$ at $\xi$.   So we can find a compact subset $K$ of $\Omega$ such that whenever $z\in \Omega\setminus K$,
\begin{itemize}
 \item $\delta_{\Omega}(z) < \min(1, t_0)$; and
 \item For any point $\xi(z)\in \partial\Omega$ that satisfies $\delta_{\Omega}(z) = d_{\Euc}(z, \xi(z))$, given
 any vector $v\neq 0$ parallel to the supporting complex hyperplane of $\Omega$ at $\xi(z)$, the complex
 line of the form $ \uni^{\xi(z)}(z + \Cb{v})$ satisfies
 \begin{align*}
  \uni^{\xi(z)}(z + \Cb{v}) = (0,\dots,0,i\bcdot \delta_{\Omega}(z)) + \Cb{\uni_{\xi(z)}(v)}
 \end{align*}
 and intersects $\boldsymbol{{\sf M}}$ in a circle of radius $\varPsi^{-1}(\delta_{\Omega}(z))$.
\end{itemize}
Here $\uni^{\xi(z)}$ is as described in Definition~\ref{def:supp}. From this point we can argue as in the proof of \cite[Lemma~3.2]{B2016}, \emph{mutatis mutandis}, to get
\begin{align*}
 r_{\Omega}(z; v) \leq 2\varPsi^{-1}(\delta_{\Omega}(z))
\end{align*}
for each $z\in \Omega\setminus K$ and $v\in \Cb^d\setminus\{0\}$ (the purpose of $t_0$ given by \eqref{eq:cap} is to ensure that
$\delta_{\Omega}(z)$ is in the domain of $\varPsi^{-1}$ and that $\delta_{\Omega}(z)\leq \varPsi^{-1}(\delta_{\Omega}(z))$ for the aforementioned $z$).

Therefore, from \eqref{eq:kob_bounds}, we deduce that 
\begin{align*}
  \frac{1}{k_{\Omega}(z; v)} \leq 4\varPsi^{-1}(\delta_{\Omega}(z))
\end{align*}
for each $z\in \Omega\setminus K$ and $v\in \Cb^d$ such that $\norm{v} = 1$. Therefore, writing $\epsilon^* := \min_{z\in K}\delta_{\Omega}(z)$, we have
\begin{align*}
 M_{\Omega}(t) \leq 4\varPsi^{-1}(t) \text{ for } t < \epsilon^*
\end{align*}
(by construction, $0<\epsilon^*\leq \epsilon$). Hence, by hypothesis, $\Omega$ satisfies Condition~1.
\end{proof}

\begin{remark}
The last lemma expresses quantitatively the claim that, for a convex domain $\Omega\Subset \Cb^d$ that satisfies Condition~1, $\partial\Omega$ is allowed to osculate $H_\xi(\partial\Omega)$ to infinite order at a smooth point $\xi\in \partial\Omega$. The function $\varPsi$ in Lemma~\ref{lem:convex_cond_1} also gives a sufficient condition on the extent to which the boundary of $\Omega$ must bend at a point $\xi\in \partial\Omega$ of infinite type for Condition~1 to hold. We illustrate all this via a familiar family of functions on $[0,\infty)$ that vanish to infinite order at $0$: these are the functions $\varPsi_s$, $s > 0$:
\begin{equation*}
 \varPsi_s(t) := \begin{cases}
 				e^{-1/t^s},	&\text{if $t > 0$}, \\
 				0,			&\text{if $t = 0$}.
 				\end{cases}
\end{equation*}
The description of the range of $s$ for which
\begin{align*}
\int_{0}^1 t^{-1}\varPsi_s^{-1}(t)\,dt = \int_0^1 t^{-1}\big(\log(1/t)\big)^{-1/s}\,dt < \infty
\end{align*}
is a standard result; $\varPsi_s$ (restricted to a suitably small interval) satisfies the conditions of the above lemma if and only if $0 < s < 1$.  
 \end{remark}

\subsection{Implications for holomorphic mappings.}\label{ssec:implications}
We shall now reformulate some of the results stated in Section~\ref{sec:intro} for the domains discussed in subsections~\ref{ssec:cond_2}~and~\ref{ssec:cond_1} and compare these results to the state of the art for certain problems of continuing interest.

Recall the works cited in subsection~\ref{ssec_WD} in connection with Wolff--Denjoy-type theorems in higher dimensions. All of the results in those works that concern the study of the iterates of a fixed-point-free holomorphic self-map on a domain, and the conclusion that this \emph{entire} sequence converges locally uniformly to a constant map, involve a domain $\Omega\Subset \Cb^d$ that:   
\begin{itemize}
 \item is topologically contractible; and
 \item is such that $\partial\Omega$ satisfies some non-degeneracy condition: some form of strict convexity in \cite{A1988, B2012, A2014}; strong pseudoconvexity in \cite{H1984}.
\end{itemize}
Call the limit point (which lies in $\partial\Omega$) appearing in all of the Wolff--Denjoy-type results just cited a \emph{Wolff point}. 
Now, it is not hard to see that an attempt to extend the dichotomy presented by the classical Wolff--Denjoy theorem (i.e., Result~\ref{res:WD} above) to higher dimensions will fail if the domain in question is not contractible. Nevertheless, it would be interesting if\,---\,making reasonable assumptions on the domain $\Omega$, but without assuming contractibility\,---\,there were to be a dichotomy wherein one of the possibilities is that the entire sequence of iterates of a holomorphic self-map of $\Omega$ converges locally uniformly to a Wolff point. In this circumstance, Theorem~\ref{thm:WD} presents the right dichotomy.

It would be even more interesting if the latter dichotomy could be exhibited for weakly pseudoconvex domains: almost none of the methods in \cite{H1984} is usable if the domain in question is a non-convex weakly pseudoconvex domain, \emph{even} if it is of finite type. A bounded pseudoconvex domain with $C^1$ boundary is taut; see~\cite{KR1981}. Thus, in view of Lemma~\ref{lem:int_cone} and the discussion in subsections~\ref{ssec:cond_2}~\&~\ref{ssec:cond_1}, Theorem~\ref{thm:WD} gives us:

\begin{corollary}\label{cor:WD_finite-type}
Let $\Omega\subset \Cb^d$ be a bounded pseudoconvex domain of finite type. If $f:\Omega\rightarrow \Omega$ is a holomorphic map then either:
\begin{enumerate}
\item for any $x \in \Omega$ the orbit $\{ f^n(x): n \in \Nb\}$ is relatively compact in $\Omega$; or
\item there exists $\xi \in \partial \Omega$ such that
\begin{equation*}
\lim_{n \rightarrow \infty} f^n(x) = \xi
\end{equation*}
 for any $x \in \Omega$, this convergence being uniform on compact subsets of $\Omega$.
\end{enumerate}
\end{corollary}

\begin{remark} Karlsson gave a proof of the above Corollary with the additional assumption that $(\Omega, K_\Omega)$ is Cauchy complete~\cite[Theorem 3]{K2005b}. This assumption greatly simplifies the situation.   \end{remark}

The discussion in subsection~\ref{ssec:cont_extn} allows us to improve upon what is currently known about the continuous extension of of a complex geodesic $\varphi : \Delta\rightarrow \Omega$ up to $\partial\Delta$, where $\Omega$ is any domain that admits complex geodesics. By Royden--Wong \cite{RW1983}, for any pair of points of a bounded convex domain $\Omega\subset \Cb^d$, $d\geq 2$\,---\,with no constraint on the regularity of $\partial\Omega$\,---\,there exists a complex geodesic of $\Omega$ containing these two points. Lempert \cite{L1984} has shown an analogous result for strongly linearly convex domains in $\Cb^d$ with $C^\infty$-smooth boundary. The result has been proved in \cite{L1984} for domains with real-analytic boundary, but the arguments therein can be adapted to the smooth case; also see \cite{KW2013}. We refer the reader to \cite{L1984} or \cite{KW2013} for a definition of strong linear convexity. It follows from the discussion  on smoothly-bounded Hartogs domains in \cite[Chapter~2]{APS2004} that strongly linearly convex domains need not necessarily be convex. Recently, Pflug and Zwonek \cite{PZ2012} provided explicit examples of strongly linearly convex domains that are not even biholomorphic to any convex domain. However, a strongly linearly convex domain is always strongly pseudoconvex; see \cite[Propositions~2.1.8 and 2.5.9]{APS2004}.

In the works cited in subsection~\ref{ssec:cont_extn} in connection with boundary regularity of complex geodesics, the domains considered were convex domains with boundaries having some degree of smoothness. Owing to Theorem~\ref{thm:qi_ext} we are able to extend those results to certain convex domains with non-smooth boundary. In \cite{L1984}, Lempert showed that in a strongly linearly convex domain with \emph{real-analytic} boundary, all complex geodesics extend real-analytically to $\partial\Delta$. However, this has a difficult and technical proof, and the proof of even $C^{1/2}$ extension is hard. The analogous result for strongly linearly convex domains with $C^\infty$ boundary is expected to have an even more technical proof. In contrast, in the $C^\infty$ setting our methods provide a rather ``soft'' proof of the continuous extension of complex geodesics up to $\partial\Delta$. To be more precise:   
owing to Theorem~\ref{thm:qi_ext}, the discussion in subsections~\ref{ssec:cond_2}~\&~\ref{ssec:cond_1}, and the fact that $(\Delta, \rho_\Delta)$ is Gromov hyperbolic, the following corollary is immediate:

\begin{corollary}
Let $\Omega\subset \Cb^d$, $d\geq 2$, be a bounded domain having either one of the following properties:
\begin{enumerate}
 \item $\Omega$ is a convex Goldilocks domain (for instance, $\Omega$ is a domain of finite type, or satisfies the conditions in Lemma~\ref{lem:convex_cond_1}); or
 \item $\Omega$ is a smoothly bounded strongly linearly convex domain.
\end{enumerate}
Then every complex geodesic $\varphi : \Delta\rightarrow \Omega$ extends to a continuous map $\wt{\varphi}: \overline{\Delta}\rightarrow \overline{\Omega}$.
\end{corollary}

\subsection{Goldilocks domains are pseudoconvex.}\label{ssec:pseudoconvex}
All the classes of domains presented above were examples of pseudoconvex domains. This is no coincidence: as hinted at in Remark~\ref{rem:just_right}, Goldilocks domains are necessarily pseudoconvex. We shall now present a proof of this. To do so, we refer to a classical result:

\begin{result}\label{res:con_principle}
A domain $\Omega\subset \Cb^d$, $d\geq 2$, is pseudoconvex if and only if for any continuous family of analytic discs\,---\,i.e., any continuous map $\Phi : \overline{\Delta}\times [0,1]\rightarrow \Cb^d$ such that $\varphi_t := \left.\Phi(\bcdot, t)\right|_{\Delta}$ is holomorphic for each $t\in [0,1]$\,---\,that satisfies $\Phi(\Delta\times\{0\}\cup \partial\Delta\times[0,1])\subset \Omega$, it follows that $\varphi_t(\Delta)\subset \Omega$ for each $t\in [0,1]$.
\end{result}

It is known, and can be ascertained by working through each step of the proof of Result~\ref{res:con_principle} that, in this result, it suffices to consider \emph{special} continuous families of analytic discs for which:
\begin{itemize}
 \item[$(a)$] each $\varphi_t$, $t\in [0,1]$, is a holomorphic immersion of $\Delta$ into $\Cb$; and
 \item[$(b)$] there exists a constant $c > 0$ such that $\|\varphi_t^\prime(\zeta)\| > c$ for every $(\zeta, t)\in \Delta\times[0,1]$.
\end{itemize}
The above can be deduced, for instance, using an intermediate characterization for pseudoconvex domains involving the so-called Hartogs figures\,---\,see, for instance, Chapter~II, \S1 of the book \cite{FG2002} by Fritzsche--Grauert. Since we do not require pseudoconvexity of Goldilocks domains in any of our proofs, we shall not elaborate on the above point any further. With this we can give a proof.

\begin{proposition}
If $\Omega \subset \Cb^d$ is a Goldilocks domain, then $\Omega$ is pseudoconvex.
\end{proposition} 

\begin{proof} Since planar domains are pseudoconvex, we consider the case $d\geq 2$.
Let $\Omega\subset \Cb^d$, $d\geq 2$, be a bounded domain. Suppose $\Omega$ is \emph{not} pseudoconvex.
There exists a continuous family of analytic discs $\Phi : \overline{\Delta}\times [0,1]\rightarrow \Cb^d$ satisfying the hypothesis of Result~\ref{res:con_principle} and conditions $(a)$ and $(b)$, but such that the conclusion of Result~\ref{res:con_principle} fails. Let
\[
 \tau\,:=\,\inf\{t\in (0,1] : \varphi_t(\Delta)\not\subset \Omega\}.
\]
As $\Omega$ is not pseudoconvex, and as the condition $\varphi_t(\overline{\Delta})\subset \Omega$ is an open condition, $\tau\in (0,1]$. By definition, there exists a point $\xi\in \partial\Omega$ and a point $\zeta_0\in \Delta$ such that $\varphi_\tau(\zeta_0) = \xi$. For $\nu\in \Zb_+$ so large that $(\tau - 1/\nu)\in (0,\tau)$, write $z_\nu := \varphi_{\tau-(1/\nu)}(\zeta_0)$. By the continuity $\Phi$,
\begin{equation}\label{eq:bdry_limit}
 z_\nu\rightarrow \xi \; \; \text{as $\nu\to +\infty$}.
\end{equation}
Let $v_\nu$ be a unit vector such that $\varphi_{\tau-(1/\nu)}^\prime(\zeta_0)\in \Cb\!v_\nu$. By the definition of the infinitesimal Kobayashi metric (owing to which it is contracted by holomorphic maps), we have
\begin{align*}
 k_{\Omega}(z_\nu; v_\nu)
 = \frac{k_{\Omega}(\varphi_{\tau-(1/\nu)}(\zeta_0);\,\varphi_{\tau-(1/\nu)}^\prime(\zeta_0))}
 	{\|\varphi_{\tau-(1/\nu)}^\prime(\zeta_0)\|} 
 &\leq \frac{k_{\Delta}(\zeta_0; 1)}{\|\varphi_{\tau-(1/\nu)}^\prime(\zeta_0)\|} \\
 &=\frac{1}{\|\varphi_{\tau-(1/\nu)}^\prime(\zeta_0)\|(1 - |\zeta_0|^2)} \\
 &\leq 1/c (1 - |\zeta_0|^2).
\end{align*}
Owing to \eqref{eq:bdry_limit} we conclude that there is an $\epsilon > 0$ such that $M_{\Omega}(r) \geq c(1-|\zeta_0|^2)$ for any $r \in (0, \epsilon)$. But this implies that Condition~1 in Definition~\ref{def:good_domains} does nor hold in $\Omega$. In particular, $\Omega$ is not a Goldilocks domain\,---\,which establishes the result.
\end{proof}

\section{Preliminaries}

\subsection{The Kobayashi distance and metric} Let $\Omega \subset \Cb^d$ be a domain. We assume that the reader is familiar with the definitions of the Kobayashi pseudo-distance $K_{\Omega}$ and the Kobayashi--Royden pseudo-metric $k_{\Omega}$ on $\Omega$. It turns out that $K_{\Omega}$ is the integrated form of $k_{\Omega}$, but this is a \emph{result} stemming from the definitions of $K_{\Omega}$ and $k_{\Omega}$; see Result~\ref{res:local_global} below. Since we shall require the original definition of $K_{\Omega}$ in a few arguments below, we give this definition. Given points $x, y\in \Omega$, we define
\begin{align*}
 K_{\Omega}(x, y) := \inf\left\{\sum_{i = 1}^n\rho_{\Delta}(\zeta_{i-1}, \zeta_i) : (\phi_1,\dots, \phi_n; \zeta_0,\dots, \zeta_n)
 \in \mathfrak{A}(x, y)\right\}
\end{align*}
where $\mathfrak{A}(x, y)$ is the set of all analytic chains in $\Omega$ joining $x$ to $y$. Here, $(\phi_1,\dots, \phi_n; \zeta_0,\dots, \zeta_n)$ is an analytic chain in $\Omega$ joining $x$ to $y$ if $\phi_i\in \Hol(\Delta, \Omega)$ for each $i$,
\begin{align*}
 x =&\,\phi_1(\zeta_0), \; \; \; \; \; \phi_n(\zeta_n) = y, \; \; \text{and} \\
 &\,\phi_i(\zeta_i) = \phi_{i+1}(\zeta_i) 
\end{align*}
for $i = 1,\dots n-1$.

Now suppose $\Omega \subset \Cb^d$ is a bounded domain. In that case,
the Kobayashi pseudo-distance is a true distance and the Kobayashi--Royden pseudo-metric is a metric. Royden~\cite[Proposition 3]{R1971} proved that the function $k_\Omega$ is upper-semicontinuous. So if a path $\sigma: [a,b] \rightarrow \Omega$ is absolutely continuous (as a map $[a,b] \rightarrow \Cb^d$) then the function $[a,b]\ni t \mapsto k_\Omega(\sigma(t); \sigma^\prime(t))$ is integrable and we can define the \emph{length of $\sigma$} to  be
\begin{align*}
\ell_\Omega(\sigma)= \int_a^b k_\Omega(\sigma(t); \sigma^\prime(t)) dt.
\end{align*}
The Kobayashi metric has the following connections to the Kobayashi distance:

\begin{result}\label{res:local_global}
 Let $\Omega \subset \Cb^d$ be a bounded domain.
\begin{enumerate}
\item \cite[Theorem 1(ii)]{NP2008} Suppose, for a point $x\in \Omega$, $k_\Omega(x;\,\bcdot)$ is continuous and positive on $\Cb^d\setminus\{0\}$. Then
\[
k_\Omega(x;v) = \lim_{h \rightarrow 0} \frac{1}{\abs{h}} K_\Omega(x,x+hv).
\]
\item \cite[Theorem 1]{R1971} For any $x,y \in \Omega$ we have
\begin{multline*}
 K_\Omega(x,y) = \inf \left\{\ell_\Omega(\sigma)\,:\,\sigma\!:\![a,b]
 \rightarrow \Omega \text{ is piecewise } C^1,\right. \\
 \left. \text{ with } \sigma(a)=x, \text{ and } \sigma(b)=y\right\}.
\end{multline*}
\item \cite[Theorem 3.1]{V1989} For any $x,y \in \Omega$ we have
\begin{multline*}
 K_\Omega(x,y) = \inf \left\{\ell_\Omega(\sigma) : \sigma\!:\![a,b]
 \rightarrow \Omega \text{ is absolutely continuous}, \right. \\
 \left. \text{ with } \sigma(a)=x, \text{ and } \sigma(b)=y\right\}.
\end{multline*}

\end{enumerate}

\end{result}

\begin{remark}
The first result above is a weaker version\,---\,which suffices for our purposes\,---\,of a result by Nikolov and Pflug \cite{NP2008}. Among other things, their result holds true on complex manifolds.
\end{remark}

\subsection{The Hopf--Rinow theorem}

Given a metric space $(X,d)$, the length of a continuous curve $\sigma:[a,b] \rightarrow X$ is defined to be
\begin{align*}
 L_d(\sigma) = \sup \left\{ \sum_{i=1}^n d(\sigma(t_{i-1}), \sigma(t_i) ): a = t_0 < t_2 < \dots < t_n=b\right\}.
\end{align*}
Then the induced metric $d_I$ on $X$ is defined to be 
\begin{align*}
d_I(x,y) = \inf\left\{ L_d(\sigma) : \sigma\!:\![a,b] \rightarrow X \text{ is continuous}, \sigma(a)=x, \text{ and } \sigma(b)=y\right\}.
\end{align*}
When $d_I = d$, the metric space $(X,d)$ is called a \emph{length metric space}. When the Kobayashi pseudo-distance is actually a distance, then the metric space $(\Omega, K_\Omega)$ is a length metric space (by construction). For such metric spaces we have the following characterization of Cauchy completeness:

\begin{result}[Hopf--Rinow] Suppose $(X,d)$ is a length metric space. Then the following are equivalent:
\begin{enumerate}
\item $(X,d)$ is a proper metric space; that is, every bounded set is relatively compact. 
\item $(X,d)$ is Cauchy complete and locally compact. 
\end{enumerate}
\end{result}

For a proof see, for instance,  Proposition~3.7 and Corollary~3.8 in Chapter~I of~\cite{BH1999}. 

When $\Omega \subset \Cb^d$ is a bounded domain the Kobayashi distance generates the standard topology on $\Omega$ and so the metric space $(\Omega, K_\Omega)$ is locally compact. In particular we obtain:

\begin{result}\label{res:hopf_rinow}Suppose $\Omega \subset \Cb^d$ is a bounded domain. Then the following are equivalent:
\begin{enumerate}
\item $\Omega, K_\Omega)$ is a proper metric space; that is, every bounded set is relatively compact. 
\item $(\Omega,K_\Omega)$ is Cauchy complete. 
\end{enumerate}
\end{result}

\subsection{Lipschitz continuity of the Kobayashi distance and metric}
We begin with a simple proposition. Since we shall re-use some aspects of the proof elsewhere in this paper, and the proof itself is short, we provide a proof.
 
\begin{proposition}\label{prop:lip}
Suppose $\Omega \subset \Cb^d$ is bounded domain.
\begin{enumerate}
\item There exists $c_1 > 0$ so that 
\begin{align*}
c_1 \norm{v} \leq k_\Omega(x;v)
\end{align*}
for all $x \in \Omega$ and $v\in\Cb^d$. In particular, 
\begin{align*}
c_1 \norm{x-y} \leq K_\Omega(x,y)
\end{align*}
for all $x,y \in \Omega$.
\item For any compact set $K \subset \Omega$ there exists $C_1=C_1(K) > 0$ so that 
\begin{align*}
k_\Omega(x;v) \leq C_1 \norm{v}
\end{align*}
for all $x \in K$ and $v \in \Cb^d$.
\item For any compact set $K \subset \Omega$ there exists $C_2=C_2(K) > 0$ so that 
\begin{align*}
K_\Omega(x,y) \leq C_2 \norm{x-y}
\end{align*}
for $x,y \in K$. 
\end{enumerate}
\end{proposition}

\begin{proof}
Fix $R > 0$ so that $\Omega$ is relatively compact in $B_R(0)$. Then 
\begin{align*}
c_1:=\inf_{x \in \overline{\Omega}, \norm{v}=1} \frac{k_{B_R(0)}(x;v)}{\norm{v}} \leq \inf_{x \in \Omega, \norm{v}=1} \frac{k_{\Omega}(x;v)}{\norm{v}}.
\end{align*}
The continuity of 
\begin{align*}
B_R(0)\times \Cb^d \ni (x,v) \mapsto k_{B_R(0)}(x;v)
\end{align*}
implies that $c_1 > 0$. Thus 
\begin{align*}
k_\Omega(x;v) \geq c_1 \norm{v}
\end{align*}
for all $x \in \Omega$ and $v \in \Cb^d$. Then, it follows from part~(2) of  Result~\ref{res:local_global} that
  \begin{align*}
K_\Omega(x,y) \geq c_1 \norm{x-y}
\end{align*}
for all $x,y \in \Omega$. This establishes part~(1). 

Next fix a compact set $K \subset \Omega$. Then there exists $r > 0$ so that $B_{2r}(x) \subset \Omega$ for all $x \in K$. Then 
\begin{equation}\label{eq:k-metric_optimal}
k_\Omega(x;v) \leq \frac{1}{2r} \norm{v} \; \; \; \forall x \in K \text{ and } \forall v \in \Cb^d.
\end{equation}
So part~(2) is true. Now, since $B_{2r}(x) \subset \Omega$ for all $x \in K$, we see that 
\begin{align*}
K_\Omega(x,y) \leq  K_{B_{2r}(x)}(x, y) \leq \frac{1}{r} \norm{x-y}
\end{align*}
when $x \in K$ and $y \in B_{r}(x)$. Now let 
\begin{align*}
M := \sup\{ K_\Omega(x,y) : x,y \in K\}.
\end{align*}
By the continuity of $K_\Omega$ we see that $M < \infty$. Then 
\begin{align*}
K_\Omega(x,y) \leq \max\left\{ \frac{1}{r}, \frac{M}{r} \right\} \norm{x-y}
\end{align*}
for all $x, y \in K$. This establishes part~(3). 
\end{proof}

\section{Length minimizing curves}\label{sec:curves}

Suppose $\Omega \subset \Cb^d$ is a bounded domain. If $I \subset \Rb$ is an interval, a map $\sigma: I \rightarrow \Omega$ is called a \emph{real geodesic} if 
\begin{align*}
K_\Omega(\sigma(s),\sigma(t)) = \abs{t-s}
\end{align*} 
for all $s,t \in I$. By Result~\ref{res:local_global}, for any two points $x, y \in \Omega$ there exists a sequence of curves joining $x$ and $y$ 
 whose lengths approach $K_\Omega(x,y)$. However, for a general bounded domain the metric space $(\Omega, K_\Omega)$ may not be Cauchy complete and 
thus there is no guarantee that this sequence of curves has a convergent subsequence. In particular, it is not clear if every two points in $\Omega$ are joined by a real geodesic. This possibility of non-existence motivates the next definition:

\begin{definition}\label{def:almost_geodesic}
Suppose $\Omega \subset \Cb^d$ is a bounded domain and $I \subset \Rb$ is an interval. For $\lambda \geq 1$ and $\kappa \geq 0$ a curve $\sigma:I \rightarrow \Omega$ is called an \emph{$(\lambda, \kappa)$-almost-geodesic} if 
\begin{enumerate} 
\item for all $s,t \in I$  
\begin{align*}
\frac{1}{\lambda} \abs{t-s} - \kappa \leq K_\Omega(\sigma(s), \sigma(t)) \leq \lambda \abs{t-s} +  \kappa;
\end{align*}
\item $\sigma$ is absolutely continuous (whence $\sigma^\prime(t)$ exists for almost every $t\in I$), and for almost every $t \in I$
\begin{align*}
k_\Omega(\sigma(t); \sigma^\prime(t)) \leq \lambda.
\end{align*}
\end{enumerate}
\end{definition}

In Proposition~\ref{prop:almost_geod_exist} below, we will show for any bounded domain $\Omega$ any two points $x,y \in \Omega$ can be joined by 
an $(1,\kappa)$-almost-geodesic. 

\begin{remark} For many domains inward pointing normal lines can be parametrized as $(1,\kappa)$-almost geodesics: for convex domains with $C^{1,\alpha}$ boundary this follows from ~\cite[Proposition 4.3]{Z2015} and for strongly pseudo-convex domains this follows from estimates in~\cite{FR1987}.
\end{remark}

\begin{proposition}\label{prop:ag_Lip}
Suppose $\Omega \subset \Cb^d$ is a bounded domain. For any $\lambda \geq 1$ there exists a $C = C(\lambda)>0$ so that any $(\lambda, \kappa)$-almost-geodesic $\sigma:I \rightarrow \Omega$ is $C$-Lipschitz (with respect to the Euclidean distance). 
\end{proposition}

\begin{proof}
By Proposition~\ref{prop:lip} there exists $c_1 > 0$ so that 
\begin{align*}
k_\Omega(x;v) \geq c_1 \norm{v}
\end{align*}
for all $x \in \Omega$ and $v \in \Cb^d$. We claim that every $(\lambda, \kappa)$-almost-geodesic is $\lambda/c_1$-Lipschitz (with respect to the Euclidean distance). 

Suppose that $\sigma:I \rightarrow \Omega$ is an $(\lambda, \kappa)$-almost-geodesic. Then for almost every $t \in I$ we have
\begin{align*}
\norm{\sigma^\prime(t)} \leq \frac{1}{c_1} k_\Omega(\sigma(t);\sigma^\prime(t)) \leq \frac{\lambda}{c_1}.
\end{align*}
Since $\sigma$ is absolutely continuous we have
\begin{align*}
\sigma(t) = \sigma(s) + \int_s^{t} \sigma^\prime(r) dr.
\end{align*}
Thus 
\begin{align*}
\norm{\sigma(t) - \sigma(s)} = \norm{ \int_s^t \sigma^\prime(r)dr } \leq  \frac{\lambda}{c_1} \abs{t-s}
\end{align*}
and $\sigma$ is $\lambda/c_1$-Lipschitz. 
\end{proof}

\begin{proposition}\label{prop:almost_geod_exist}
Suppose $\Omega \subset \Cb^d$ is a bounded domain. For any $\kappa > 0$ and $x,y \in \Omega$ there exists an $(1,\kappa)$-almost-geodesic $\sigma:[a,b] \rightarrow \Omega$ with $\sigma(a) =x$ and $\sigma(b)=y$. 
\end{proposition}

We begin the proof with a simple lemma:

\begin{lemma}\label{lem:restricted_l}
Suppose $\Omega \subset \Cb^d$ is a bounded domain and $\sigma: [a,b] \rightarrow \Omega$ is an absolutely continuous curve. If 
\begin{align*}
\ell_\Omega(\sigma) \leq K_\Omega(\sigma(a),\sigma(b)) + \epsilon
\end{align*}
then whenever $a \leq a^\prime \leq b^\prime \leq b$ we have
\begin{align*}
\ell_\Omega(\sigma|_{[a^\prime, b^\prime]}) \leq K_\Omega(\sigma(a^\prime),\sigma(b^\prime)) + \epsilon.
\end{align*}
\end{lemma}

\noindent{This lemma is an immediate consequence of the fact that 
\begin{align*}
 \ell_\Omega(\sigma|_{[a^\prime, b^\prime]}) 
&= \ell_\Omega(\sigma) - \ell_\Omega(\sigma|_{[a, a^\prime]}) - \ell_\Omega(\sigma|_{[b^\prime,b]}),
\end{align*}
and of the triangle inequality for $K_\Omega$.}

\begin{proof}[The proof of Proposition~\ref{prop:almost_geod_exist}]
By part ~(2) of Result~\ref{res:local_global} there exists a piecewise $C^1$ curve $\sigma:[0,T] \rightarrow \Omega$ so that $\sigma(0) = x$, $\sigma(T) = y$, and 
\begin{align*}
\ell_\Omega(\sigma) < K_\Omega(x,y)+\kappa.
\end{align*}
Since $k_\Omega$ is upper semi-continuous, we can perturb $\sigma$ and assume, in addition, that $\sigma$ is $C^1$-smooth, and that $\sigma^\prime(t) \neq 0$ for all $t \in [0,T]$. 

Next consider the function 
\begin{equation*}
f(t) = \int_0^t k_\Omega(\sigma(s); \sigma^\prime(s)) ds.
\end{equation*}
Since $\sigma([0,T])$ is compact, by Proposition~\ref{prop:lip} there exists $C \geq 1$ so that 
\begin{equation}\label{eq:deriv_bounds}
\frac{1}{C} \norm{\sigma^\prime(t)} \leq k_\Omega(\sigma(t), \sigma^\prime(t)) \leq C \norm{\sigma^\prime(t)}
\; \; \; \text{for all } t \in [0,T].
\end{equation}
 Thus, since $\sigma^\prime(t) \neq 0$ for all $t \in [0,T]$, we see that $f$ is a bi-Lipschitz strictly increasing function.

Next let $g:[0,\ell_\Omega(\sigma)] \rightarrow [0,T]$ be the inverse of $f$, that is $f(g(t)) = t$ for all $t \in [0,T]$. We claim that the curve $\sigma_0 := \sigma\circ g$ is an $(1,\kappa)$-almost-geodesic. Since $g$ is Lipschitz ($f$ is bi-Lipschitz) we see that $\sigma_0$ is Lipschitz and hence absolutely continuous. Moreover, if $g^\prime(t)$ exists then 
\begin{align*}
\sigma_0^\prime(t) = \sigma^\prime(g(t)) g^\prime(t).
\end{align*}
When $g^\prime(t)$ exists and $f^\prime(g(t))$ exists and is non-zero, we have 
\begin{align*}
g^\prime(t) = \frac{1}{f^\prime(g(t))}.
\end{align*}
Now, by the Lebesgue differentiation theorem applied to $f$, there exists a set $E \subset [0,T]$ of full measure so that if $s \in E$ then $f^\prime(s)$ exists and 
\begin{align*}
f^\prime(s) = k_\Omega(\sigma(s); \sigma^\prime(s)).
\end{align*}
Since $g$ is bi-Lipschitz, $g^{-1}(E) \subset [0,\ell_\Omega(\sigma)]$ has full measure. Hence, as $\sigma^\prime(t)\neq 0$ for all $t \in [0,T]$, we can write (in view of \eqref{eq:deriv_bounds} above)
\begin{align*}
g^\prime(t) = \frac{1}{k_\Omega(\sigma(g(t)); \sigma^\prime(g(t)))}
\end{align*}
for almost every $t \in [0,\ell_\Omega(\sigma)]$. So for almost every $t \in [0, l_\Omega(\sigma)]$
\begin{align*}
k_\Omega( \sigma_0(t); \sigma_0^\prime(t)) = k_\Omega\Big(\sigma(g(t)); \sigma^\prime(g(t))g^\prime(t)\Big)=1.
\end{align*}
Therefore
\begin{align*}
\ell_{\Omega}(\sigma_0) = \ell_\Omega(\sigma) \leq K_\Omega(x,y) + \kappa.
\end{align*}
So, by Lemma~\ref{lem:restricted_l}, whenever $0 \leq s \leq t \leq \ell_\Omega(\sigma)$ we have 
\begin{align*}
\abs{t-s} = \ell_\Omega(\sigma_0|_{[s,t]}) \leq K_\Omega(\sigma_0(t),\sigma_0(s)) +\kappa.
\end{align*}
Since $\sigma_0$ is absolutely continuous, Result~\ref{res:local_global} implies that
\begin{align*}
K_\Omega(\sigma_0(t),\sigma_0(s)) \leq  \ell_\Omega(\sigma_0|_{[s,t]}) = \abs{t-s}.
\end{align*}
So $\sigma_0$ is an $(1,\kappa)$-almost-geodesic.
\end{proof}

\subsection{Real geodesics}

In this subsection we show that when $\Omega$ is a taut bounded domain, then any real geodesic must possess a certain degree of extra regularity: namely, that it is an $(1,0)$-almost-geodesic.

\begin{proposition}\label{prop:geodesics}
Suppose $\Omega \subset \Cb^d$ is a bounded domain. Then there exists $C_\Omega > 0$ so that any real geodesic $\sigma : I \rightarrow \Omega$ is $C_\Omega$-Lipschitz  (with respect to the Euclidean distance). In particular 
\begin{equation}\label{eq:ftCalc}
\sigma(t) = \sigma(t_0) + \int_{t_0}^t \sigma^\prime(s) ds 
\end{equation}
for any $t,t_0 \in I$. Moreover, if $\Omega$ is taut then
\begin{align*}
K_\Omega(\sigma(t); \sigma^\prime(t)) = 1
\end{align*}
for almost every $t \in I$. 
\end{proposition}

\begin{proof}
By Proposition~\ref{prop:lip} there exists $c >0$ so that 
\begin{align*}
K_\Omega(x,y) \geq c \norm{x-y}
\end{align*}
for all $x,y \in \Omega$. Then
\begin{align*}
\norm{\sigma(t) - \sigma(s)} \leq \frac{1}{c} \abs{t-s}.
\end{align*}
Thus $\sigma$ is Lipschitz (and $1/c$ is the $C_{\Omega}$ mentioned above). In particular, $\sigma$ is absolutely continuous, from which \eqref{eq:ftCalc} follows.
 
Next suppose that $\Omega$ is taut. We now appeal to Theorem~1.2 in \cite{V1989}. Since $\Omega$ is a taut and bounded domain, $k_\Omega$ is continuous; see, for instance, \cite[Section~3.5]{JP1993}. Hence the conditions stated in part~(1) of Result~\ref{res:local_global} are satisfied. Thus, in our specific context, \cite[Theorem 1.2]{V1989} reads as
\begin{align*}
\int_{t}^{t+h} k_\Omega(\sigma(s); \sigma^\prime(s))ds
 = \sup_{\mathcal{P}}\sum_{j=1}^{N(\mathcal{P})} K_\Omega(\sigma(s_{j-1}), \sigma(s_j))
\end{align*}
where the supremum above ranges over all partitions
\begin{align*}
\mathcal{P}\,:\,t=s_0 < s_1 < s_2 < \dots < s_{N(\mathcal{P})}=t+h
\end{align*}
of $[t, t+h]$, and where $h > 0$ is such that $t, t+h \in I$. As $\sigma$ is a real geodesic, we then have
\begin{align*}
 \int_{t}^{t+h} k_\Omega(\sigma(s); \sigma^\prime(s))ds = h
\end{align*}
for all $h > 0$ such that $t, t+h \in I$. Then by the Lebesque differentiation theorem
\begin{align*}
k_\Omega(\sigma(t); \sigma^\prime(t)) = 1
\end{align*}
for almost every $t \in I$.

\end{proof}

\subsection{Quasi-geodesics}

In this subsection, we show that any quasi-geodesic can, in a certain sense, be approximated by an almost-geodesic. This proposition will be needed in our proof of continuous extension of isometries.

\begin{definition} Suppose $\Omega \subset \Cb^d$ is a bounded domain and $I \subset \Rb$ is a interval. For $\lambda \geq 1$, $\kappa \geq 0$ a map $\sigma:I \rightarrow \Omega$ is called a \emph{$(\lambda, \kappa)$-quasi-geodesic} if  
\begin{align*}
\frac{1}{\lambda} \abs{t-s} -\kappa \leq K_\Omega(\sigma(s), \sigma(t)) \leq \lambda \abs{t-s} +  \kappa
\end{align*}
 for all $s,t \in I$. 
\end{definition}

\begin{remark}
Note that a $(\lambda, \kappa)$-quasi-geodesic is not required to be continuous. It is in this sense that it differs from an $(\lambda, \kappa)$-\emph{almost}-geodesic, which must have greater regularity; see Definition~\ref{def:almost_geodesic}. Furthermore, we ought to remark that the proposition below makes no assertions about existence of quasi-geodesics. Also, while Proposition~\ref{prop:almost_geod_exist} asserts that any pair of points of a bounded domain $\Omega\subset \Cb^d$ can be joined by an $(1, \kappa)$-almost-geodesic\,---\,which is more regular than a $(1, \kappa)$-quasi-geodesic\,---\,this comes \emph{with the proviso that $\kappa > 0$}.
\end{remark}

\begin{proposition}\label{prop:approx_quasi_geod}
Suppose $\Omega \subset \Cb^d$ is a bounded domain. For any $\lambda \geq 1$, $\kappa \geq 0$ there exist constants
$R > 0$, $\lambda_0\geq 1$, and $\kappa_0\geq 0$, depending \emph{only} on the pair $(\lambda, \kappa)$,
that have the following property: for any $(\lambda, \kappa)$-quasi-geodesic $\sigma:[a,b] \rightarrow \Omega$ there exists an $(\lambda_0, \kappa_0)$-almost-geodesic $S : [0,T] \rightarrow \Omega$ with $S(0) = \sigma(a)$, $S(T) = \sigma(b)$, and such that
\begin{align*}
\max \left\{ \sup_{t \in [a,b]} K_\Omega(\sigma(t), S), \sup_{t \in [0,T]} K_\Omega(S(t), \sigma) \right\} \leq R.
\end{align*}
Here, given a set $E\subset \Omega$ and a point $o\in \Omega$, we write $K_{\Omega}(o, E) := \inf_{x\in E}K_{\Omega}(o, x)$.
\end{proposition}

\begin{proof}
The argument falls into two cases, depending on the magnitude of $\abs{b - a}$. The case that requires some work is when $\abs{b - a}$ is large.
\medskip

\noindent{{\bf Case 1:} First consider the case when $\abs{b - a} > 1/2$. Fix a partition
\begin{align*}
a=t_0 < t_1 < t_2 <\dots < t_{N} = b
\end{align*}
so that $1/2\leq \abs{t_{k} - t_{k-1}} \leq 1$. For $1 \leq k \leq N$, let $\gamma_k : [0,T_k] \rightarrow \Omega$ be an $(1,1)$-almost-geodesic with $\gamma_k(0) = \sigma(t_{k-1})$ and $\gamma_k(T_k)=\sigma(t_k)$; notice that such a curve exists by Proposition~\ref{prop:almost_geod_exist}.}

Now, by the properties of $\gamma_k$ we get 
\begin{align*}
T_k -1 \leq K_\Omega(\gamma_k(0), \gamma_k(T_k)) = K_\Omega(\sigma(t_{k-1}), \sigma(t_k)) \leq \lambda \abs{t_k-t_{k-1}} +\kappa \leq \lambda + \kappa,
\end{align*}
whence $T_k \leq \lambda + \kappa+1$, $k = 1,\dots, N$. Therefore, 
\begin{align*}
K_\Omega(\gamma_k(t), \sigma(t_{k-1}))
= K_\Omega(\gamma_k(t), \gamma_k(0)) &\leq \abs{t} +1 \\
& \leq T_k + 1 \leq \lambda + \kappa+2
\end{align*}
for any $t \in [0,T_k]$. 

Let $S : [a,b] \rightarrow \Omega$ be the curve defined as follows: 
\begin{align*}
S(t) = \gamma_k\left( \frac{T_k}{t_k-t_{k-1}}( t-t_{k-1})\right), \; \text{ if } \; t_{k-1} \leq t \leq t_k,  \; k = 1,\dots, N.
\end{align*}
Then, using the estimate above, for  $t \in [t_{k-1},t_k]$ we have
\begin{align*}
K_\Omega(S(t), \sigma(t)) 
& \leq K_\Omega(S(t), \sigma(t_{k-1})) + K_\Omega(\sigma(t_{k-1}), \sigma(t))  \\
& \leq \lambda +  \kappa +2+ \lambda\abs{t_{k-1}-t} +  \kappa \\
&\leq  2 \lambda + 2 \kappa+2.
\end{align*}
Write $R := 2 \lambda + 2 \kappa+2$. Then 
\begin{align*}
\abs{K_\Omega(S(t), S(s)) - K_\Omega(\sigma(t), \sigma(s))} \leq K_\Omega(S(t), \sigma(t))
+ K_\Omega(S(s), \sigma(s)) \leq 2R 
\end{align*}
and so 
\begin{align*}
\frac{1}{\lambda}\abs{t-s} - \kappa - 2R \leq K_\Omega(S(t), S(s))  \leq \lambda\abs{t-s} +\kappa + 2R.
\end{align*}
Finally since each $\gamma_k$ is an $(1,1)$-almost-geodesic we see that 
\begin{align*}
k_\Omega(S(t); S^\prime(t))\, \leq \max_{1 \leq k \leq N}\,\frac{T_k}{t_{k+1}-t_k} \leq
2\lambda + 2\kappa + 2
\end{align*}
for almost every $t \in [a,b]$. 

Thus $S:[a,b] \rightarrow \Omega$ is an $(\lambda_0, \kappa_0)$-almost-geodesic where $\lambda_0 = 2\lambda+2\kappa+2$ and $\kappa_0 = \kappa+2R = 4\lambda+5\kappa+4$.
\medskip

\noindent{{\bf Case 2:} Now consider the case when $\abs{b - a}\leq 1/2$. Let $S : [0, T]\rightarrow \Omega$ be an $(1,1)$-almost-geodesic with $S(0) = \sigma(a)$ and $S(T)=\sigma(b)$. Arguing as before shows that 
\begin{align*}
T \leq \frac{\lambda}{2} +\kappa+2.
\end{align*}
Now if $t \in [a,b]$ then
\begin{align*}
K_\Omega(\sigma(t), S(0))= K_\Omega(\sigma(t), \sigma(0)) \leq \frac{\lambda}{2} +\kappa
\end{align*}
and if $t \in [0,T]$ then 
\begin{align*}
K_\Omega(S(t), \sigma(0))= K_\Omega(S(t), S(0)) \leq T +1 \leq  \frac{\lambda}{2} +\kappa+3.
\end{align*}}
\end{proof}

\section{A visibility condition}\label{sec:visible}

This section is dedicated to proving Theorem~\ref{thm:visible}. It is a key part of the present work. What makes Theorem~\ref{thm:visible} a key part of the proofs in the later sections is that, if $\Omega$ is a Goldilocks domain, $(\Omega, K_{\Omega})$ resembles \emph{adequately} a visibility space (in the sense of \cite{EO1973}, for instance) even though $(\Omega, K_{\Omega})$ is \emph{not} in general Gromov hyperbolic, nor is it known whether every pair of points can be joined by a geodesic.

We will need the following simple observation:

\begin{observation}
Suppose $f:\Rb_{\geq 0} \rightarrow \Rb_{\geq 0}$ is a bounded Lebesgue-measurable function such that
\begin{align*}
\int\nolimits_0^{\epsilon} \frac{1}{r} f(r) dr < \infty
\end{align*}
for some (and hence any) $\epsilon >0$. Then 
\begin{align*}
\int\nolimits_0^{\infty} f(Ae^{-Bt}) dt < \infty
\end{align*}
for any $A,B > 0$. 
\end{observation}
This  is an immediate consequence of the change-of-variable formula, writing
$r = Ae^{-Bt}$ in the second integral above.

\begin{proof}[The proof of Theorem~\ref{thm:visible}] Suppose that there does not exist a compact set with the desired property. Then we can find a sequence $\sigma_n:[a_n,b_n] \rightarrow \Omega$ of $(\lambda, \kappa)$-almost-geodesics so that $\sigma_n(a_n) \in V_\xi$, $\sigma_n(b_n) \in V_\eta$, and
\begin{align*}
0 = \lim_{n \rightarrow \infty} \max\{ \delta_\Omega(\sigma_n(t)) : t \in [a_n,b_n]\}.
\end{align*}

By reparametrizing each $\sigma_n$ we can assume that 
\begin{align*}
\delta_\Omega(\sigma_n(0)) = \max\{ \delta_\Omega(\sigma_n(t)) : t \in [a_n,b_n]\}.
\end{align*}
Then by passing to a subsequence we can assume $a_n \rightarrow a \in [-\infty,0]$, $b_n \rightarrow b \in [0,\infty]$, $\sigma_n(a_n) \rightarrow \xi^\prime$, and $\sigma_n(b_n) \rightarrow \eta^\prime$. By assumption $\xi^\prime \in \overline{V_\xi} \cap \partial \Omega$ and $\eta^\prime \in \overline{V_\eta} \cap \partial \Omega$. Notice that $\xi^\prime \neq \eta^\prime$ because $\overline{V_\xi} \cap \overline{V_\eta} = \emptyset$.

By Proposition~\ref{prop:ag_Lip}, there exists some $C>0$ so that each $\sigma_n$ is $C$-Lipschitz with respect to the Euclidean distance. Thus we can pass to another subsequence so that $\sigma_n$ converges locally uniformly on $(a,b)$ to a  curve $\sigma:(a,b) \rightarrow \overline{\Omega}$ (we restrict to the open interval because $a$ could be $-\infty$ and $b$ could be $\infty$). Notice that $a \neq b$ because each $\sigma_n$ is $C$-Lipschitz and so
\begin{align*}
0 < \norm{\xi^\prime-\eta^\prime} \leq C \abs{b-a}.
\end{align*}
 
Since $\sigma_n$ is an $(\lambda, \kappa)$-almost-geodesic
\begin{equation*}
k_\Omega(\sigma_n(t); \sigma_n^\prime(t)) \leq \lambda
\end{equation*}
for almost every $t \in [a_n,b_n]$. We claim that:
\begin{equation}\label{eq:speed_estimate}
\norm{\sigma_n^\prime(t)} \leq \lambda M_\Omega(\delta_\Omega(\sigma_n(t)))
\; \; \; \text{for almost every } t \in [a_n,b_n]. 
\end{equation}
In the case when $\sigma_n^\prime(t) = 0$ this is immediate and if $\sigma_n^\prime(t) \neq 0$ we have
\begin{align*}
\norm{\sigma_n^\prime(t)} \leq \frac{\lambda}{k_\Omega\left(\sigma_n(t); \frac{1}{\norm{\sigma_n^\prime(t)}} \sigma^\prime(t) \right)} \leq \lambda M_\Omega(\delta_\Omega(\sigma_n(t))).
\end{align*}
\smallskip

\noindent \textbf{Claim 1:} $\sigma:(a,b) \rightarrow \overline{\Omega}$ is a constant map.

\noindent \emph{Proof.} Since 
\begin{align*}
\delta_\Omega(\sigma_n(t)) \leq \delta_\Omega(\sigma_n(0))
\end{align*}
we see that 
\begin{align*}
M_\Omega(\delta_\Omega(\sigma_n(t))) \leq M_\Omega(\delta_\Omega(\sigma_n(0))).
\end{align*}
Thus $M_\Omega(\delta_\Omega(\sigma_n(t))) \rightarrow 0$ uniformly. But then  if $u \leq w$ and $u,w \in (a,b)$
\begin{align*}
 \norm{\sigma(u)-\sigma(w)} 
 = \lim_{n \rightarrow \infty} \norm{\sigma_n(u)-\sigma_n(w)}
 &\leq \limsup_{n \rightarrow \infty} \int_u^w \norm{\sigma_n^\prime(t)} dt\\
 & \leq  \lambda\limsup_{n \rightarrow \infty} \int_u^w\!\!M_\Omega(\delta_\Omega(\sigma_n(t))) dt = 0.
\end{align*}
Thus $\sigma$ is constant. \hfill $\blacktriangleleft$
\medskip

We will establish a contradiction by proving the following:
\medskip

\noindent \textbf{Claim 2:} $\sigma:(a,b) \rightarrow \overline{\Omega}$ is not a constant map. 

\noindent \emph{Proof.} Fix $x_0 \in \Omega$. Then there exists $C,\alpha > 0$ so that 
\begin{align*}
K_\Omega(x,x_0) \leq C + \alpha \log \frac{1}{\delta_\Omega(z)}
\end{align*}
for all $x \in \Omega$. Therefore
\begin{align*}
\frac{1}{\lambda}\abs{t} - \kappa \leq K_\Omega(\sigma_n(0), \sigma_n(t))
& \leq K_\Omega(\sigma_n(0),x_0) + K_\Omega(x_0, \sigma_n(t)) \\
& \leq 2C + \alpha \log \frac{1}{\delta_\Omega(\sigma_n(0))\delta_\Omega(\sigma_n(t))}.
\end{align*}
Thus 
\begin{align*}
\delta_\Omega(\sigma_n(t)) \leq \sqrt{ \delta_\Omega(\sigma_n(0)) \delta_\Omega(\sigma_n(t)) } \leq A e^{-B\abs{t}}
\end{align*}
where $A = e^{(2C+\kappa)/(2\alpha)}$ and $B = 1/(2\alpha\lambda)$. 

Thus, by the estimate \eqref{eq:speed_estimate}, for almost every $t \in [a_n,b_n]$ we have
\begin{align*}
\norm{\sigma_n^\prime(t)} \leq \lambda M_\Omega(\delta_\Omega(\sigma_n(t))) \leq \lambda M_\Omega( A e^{-B\abs{ t}})
\end{align*}

Now fix $a^\prime, b^\prime \in (a,b)$ so that 
\begin{align*}
\norm{\xi^\prime-\eta^\prime} > \lambda \int_a^{a^\prime} M_\Omega(Ae^{-B\abs{t}}) dt +\lambda \int_{b^\prime}^b M_\Omega(Ae^{-B\abs{t}}) dt.
\end{align*}
Then 
\begin{align*}
\norm{\sigma(b^\prime)-\sigma(a^\prime) } 
&= \lim_{n \rightarrow \infty} \norm{\sigma_n(b^\prime)-\sigma_n(a^\prime) } \\
& \geq \lim_{n \rightarrow \infty} \big(\norm{\sigma_n(b_n)-\sigma_n(a_n) } - 
\norm{\sigma_n(b_n)-\sigma_n(b^\prime) } - \norm{\sigma_n(a^\prime)-\sigma_n(a_n) }\big) \\
& \geq \norm{\xi^\prime-\eta^\prime} - \limsup_{n \rightarrow \infty} \int_{b^\prime}^{b_n} \norm{\sigma_n^\prime(t)} dt - \limsup_{n \rightarrow \infty} \int_a^{a^\prime} \norm{\sigma_n^\prime(t)} dt \\
& \geq  \norm{\xi^\prime-\eta^\prime} - \limsup_{n \rightarrow \infty}\lambda \int_{b^\prime}^{b_n} M_\Omega(Ae^{-B \abs{t}})dt - \limsup_{n \rightarrow \infty} \lambda\int_{a_n}^{a^\prime}M_\Omega(Ae^{-B \abs{t}}) dt\\
& = \norm{\xi^\prime-\eta^\prime} - \lambda\int_{b^\prime}^b M_\Omega(Ae^{-B\abs{t}}) dt -\lambda \int_a^{a^\prime} M_\Omega(Ae^{-B\abs{t}}) dt >0.
\end{align*}
Thus $\sigma:(a,b) \rightarrow \overline{\Omega}$ is non-constant. \hfill $\blacktriangleleft$
\medskip

The above contradicts Claim~1. This establishes the existence of the compact $K$ with the stated property.
\end{proof}

\section{Extensions of quasi-isometries}\label{sec:gromov_prod}

\subsection{The Gromov boundary} Let $(X,d)$ be a metric space. Given three points $x,y,o \in X$, the \emph{Gromov product} is given by 
\begin{align*}
 (x|y)_o = \frac{1}{2} \left( d(x,o)+d(o,y)-d(x,y) \right).
\end{align*}
When $(X,d)$ is proper and Gromov hyperbolic, the Gromov product can be used to define an abstract boundary at infinity denoted $X(\infty)$ and called the \emph{Gromov boundary}. In particular, a seqeunce $(x_n)_{n \in \Nb} \subset X$ is said to \emph{converge at $\infty$} if 
\begin{align*}
\liminf_{n,m \rightarrow \infty} (x_n|x_m)_o = \infty
\end{align*}
for some (and hence) any $o \in X$. Two sequences $(x_n)_{n \in \Nb}$ and $(y_n)_{n \in \Nb}$ in $X$ are equivalent if 
\begin{align*}
\liminf_{n,m \rightarrow \infty} (x_n|y_m)_o = \infty
\end{align*}
for some (and hence) any $o \in X$. Finally $X(\infty)$ is the set of equivalence classes of sequences converging to infinity. Moreover, $X \cup X(\infty)$ has a natural topology (see for instance~\cite[Part III.H.3]{BH1999}) that makes it a compactification of $X$.

\subsection{Continuous extensions of quasi-isometries}\label{ssec:all_about_qi}

Given a bounded domain $\Omega\subset \Cb^d$, it is, in general, very hard to determine whether $(\Omega, K_{\Omega})$ is Gromov hyperbolic. Furthermore, as we saw in subsection~\ref{ssec:cont_extn}, $K_{\Omega}$ fails to be Gromov hyperbolic for domains as regular as convex domains with $C^\infty$-smooth boundary if $\partial\Omega$ has points of  infinite type; see \cite{Z2014}. This renders unusable a very natural approach, namely Result~\ref{res:gromov_qi_ext}, for studying the boundary behavior of continuous quasi-isometries (for the Kobayashi metric), even if they are holomorphic, on such domains. We therefore explore alternative notions of good compactifications\,---\,from the viewpoint of obtaining continuous extensions of quasi-isometries\,---\,of $(\Omega, K_{\Omega})$. To this end, we begin with a couple of very general definitions. 

\begin{definition}
Let $(X,d)$ be a metric space. A pair $(\iota, X^*)$ is a \emph{compactification} of $X$ if $X^*$ is a sequentially compact Hausdorff topological space, $\iota:X \rightarrow X^*$ is a homeomorphism onto its image, and $\iota(X)$ is open and dense in $X^*$.
 \end{definition}
 
\begin{definition}\label{def:gooc}
Suppose $(\iota,X^*)$ is a compactification of a geodesic metric space $(X,d)$. We say $(\iota,X^*)$ is a \emph{good compactification} if for all sequences $\sigma_n:[a_n,b_n] \rightarrow X$ of geodesics with the property
\begin{align*}
\lim_{n \rightarrow \infty} \iota(\sigma_n(a_n)) = \lim_{m \rightarrow \infty} \iota(\sigma_n(b_n)) \in X^* \setminus \iota(X)
\end{align*}
we have
\begin{align*}
\liminf_{n \rightarrow \infty} d(o, \sigma_n) = \infty
\end{align*}
for any $o \in X$.
\end{definition}

To clarify our notation: if $\sigma : [0, T]\rightarrow X$ is a map and $o\in X$,
$d(o,\sigma) := \inf\{d(o,\sigma(s)) : s \in [0,T]\}$.

As the next observation shows, good compactifications only exist for proper metric spaces:

\begin{observation} Suppose $(\iota,X^*)$ is a good compactification of a metric space $(X,d)$. Then $(X,d)$ is a proper metric space. 
\end{observation}

\begin{proof}
Fix $R > 0$ and $x_0 \in X$. We claim that the set $B:=\{ y \in X : d(y,x_0) \leq R\}$ is compact. To see this, fix a sequence $x_n \in B$. Since $X^*$ is sequentially compact we can assume, passing to a subsequence if necessary, that $\iota(x_n) \rightarrow \xi \in X^*$. If $\xi \in X$, then $\xi \in B$. If $\xi \in X^* \setminus \iota(X)$ then the curve $\sigma_n:[0,0] \rightarrow X$ given by $\sigma_n(0)=x_n$ is a geodesic. So, by the definition of a good compactification, 
\begin{align*}
\infty = \liminf_{n \rightarrow \infty} d(x_0, \sigma_n) = \liminf_{n \rightarrow \infty} d(x_0, x_n) \leq R,
\end{align*}
which is a contradiction, whence $\xi\notin X^*\setminus \iota(X)$. 
\end{proof}

\begin{remark}
We now discuss a few examples and look at some motivations underlying Definition~\ref{def:gooc}.
\begin{enumerate}
\item Let $X^* = \Rb^d \cup \{\infty\}$ be the one point compactification of $(\Rb^d, d_{\Euc})$, then $X^*$ is not a good compactification.  
\item In view of Theorem~\ref{thm:quasi_isometry_ext} below, it would be desirable if the Gromov compactification $X\cup X(\infty)$, where $(X,d)$ is a proper geodesic Gromov hyperbolic space, were subsumed by Definition~\ref{def:gooc}. This is in fact the case by Result~\ref{res:gromov_qi_ext}.
\item Let $\Omega$ be a bounded convex domain with $C^{1,\alpha}$-smooth boundary and assume that for each $\xi\in \partial\Omega$, the affine set $\xi+H_\xi(\partial\Omega)$ (see subsection~\ref{ssec:cond_1} for the definition) intersects
$\overline\Omega$ precisely at $\xi$. It is a classical fact that $(\Omega, K_{\Omega})$ is Cauchy complete, see for instance~\cite[Proposition 2.3.45]{A1989}. It then follows that  $(\Omega, K_{\Omega})$ is a geodesic metric space. If $\partial\Omega$ contains points that are not of finite type (in the sense of D'Angelo), then $(\Omega, K_{\Omega})$ is not Gromov hyperbolic; see \cite{Z2014}. Yet, irrespective of whether or not $(\Omega, K_{\Omega})$ is Gromov hyperbolic, it follows from \cite[Theorem~2.11]{Z2015} that $({\sf id}_{\Omega}, \overline\Omega)$ is a good compactification.
\end{enumerate}
\end{remark}

The next theorem could be stated for any geodesic metric space $(X, d)$ that admits a good compactification $(\iota, X^*)$ and any quasi-isometric embedding $F : (X, d) \rightarrow (\Omega, K_\Omega)$. However, it is unclear what the interest in such a general set-up could be. On the other hand, we have seen in the discussion in subsection~\ref{ssec:cont_extn} that quasi-isometric embeddings\,---\,relative to the Kobayashi metric\,---\,between \emph{domains} arise rather naturally, while existing tools for studying their boundary are no longer effective. These are the considerations behind the statement about quasi-isometries between two domains in Theorem~\ref{thm:quasi_isometry_ext}. Observe that Theorem~\ref{thm:qi_ext} is a special case of Theorem~\ref{thm:quasi_isometry_ext}.

\begin{theorem}\label{thm:quasi_isometry_ext}
Let $D$ be a bounded domain in $\Cb^k$ and suppose $(D, K_D)$ admits a good compactification $(\iota, D^*)$. Let $\Omega\subset \Cb^d$ be a Goldilocks domain. If $F : (D, K_D) \rightarrow (\Omega, K_\Omega)$ is a continuous quasi-isometric embedding, then $F$ extends to a continuous map from $D^*$ to $\overline{\Omega}$.
\end{theorem}

The following lemma is the key to proving Theorem~\ref{thm:quasi_isometry_ext}. Its proof follows immediately from  
Theorem~\ref{thm:visible} and Proposition~\ref{prop:approx_quasi_geod}.

\begin{proposition}\label{prop:qg_visibility}
 Suppose $\Omega \subset \Cb^d$ is a Goldilocks domain and $\lambda \geq1$, $\kappa \geq 0$. If $\xi,\eta \in \partial\Omega$ and $V_\xi, V_\eta$ are neighborhoods of $\xi,\eta$ in $\overline{\Omega}$ so that $\overline{V_\xi} \cap \overline{V_\eta} = \emptyset$, then for each $x_0 \in \Omega$ there exists $R > 0$ with the following property: if $\sigma: [a,b] \rightarrow \Omega$ is a $(\lambda, \kappa)$-quasi-geodesic with $\sigma(a) \in V_\xi$ and $\sigma(b) \in V_\eta$  then 
 \begin{align*}
 K_\Omega( x_0, \sigma) \leq R.
 \end{align*}
\end{proposition}

\begin{remark} If $(\Omega, K_\Omega)$ is Cauchy complete then the conclusion of Theorem~\ref{thm:visible} and Proposition~\ref{prop:qg_visibility} are equivalent (by Result~\ref{res:hopf_rinow}). But in general, Proposition~\ref{prop:qg_visibility} is weaker. This is due to the following hypothetical example:  suppose there exists two sequences $x_n \rightarrow \xi \in \partial \Omega$ and $y_n \rightarrow \eta \in \partial \Omega$ so that 
\begin{align*}
\sup_{n \in \Nb} K_\Omega(x_n, x_0) = R_1 < \infty \text{ and } \sup_{n \in \Nb}  K_\Omega(y_n, x_0) = R_2 < \infty.
\end{align*}
Then the sequence of maps $\sigma_n :[0,1] \rightarrow \Omega$ given by
\begin{align*}
\sigma_n(t) = \left\{ \begin{array}{ll} x_n & \text{ if } 0 \leq t < 1/2 \\
y_n & \text{ if } 1/2 \leq t \leq 1.
\end{array}
\right.
\end{align*}
are all $(1,R_1+R_2+1)$-quasi-geodesics. But 
\begin{align*}
\lim_{n \rightarrow \infty} \left(\max_{t \in [a_n,b_n]} \delta_\Omega(\sigma_n(t))\right)=0.
\end{align*}
\end{remark}

Theorem~\ref{thm:quasi_isometry_ext} is an application of Theorem~\ref{thm:visible}, with Proposition~\ref{prop:qg_visibility} serving as a visibility theorem for quasi-geodesics. In fact, visibility may be seen as a tool for controlling the oscillation of $F$ ($F$ as in Theorem~\ref{thm:quasi_isometry_ext}) along various sequences $(x_n)_{n \in \Nb} \subset D$ as $\iota(x_n)$ approaches some chosen point in $\xi\in D^*\setminus \iota(D)$. The idea of the proof is as follows. If the cluster set of values as one approaches $\xi$ were non-trivial, we would have two sequences $(x_n)_{n \in \Nb}$ and $(y_n)_{n \in \Nb}$ as above such that $(F(x_n))_{n \in \Nb}$ and $(F(y_n))_{n \in \Nb}$ approach two \emph{different} points in $\partial\Omega$. Let $\sigma_n$ be a geodesic joining $x_n$ to $y_n$. Then $F\circ\sigma_n$ are quasi-geodesics, whence, by Proposition~\ref{prop:qg_visibility}, these curves must be within some finite Kobayashi distance from any chosen point in $\Omega$. But then the curves $\sigma_n$ would have the analogous property in $D$, which is ruled out by the geometry of $D$.

\begin{proof}[The proof of Theorem~\ref{thm:quasi_isometry_ext}]
Fix some $\xi \in D^* \setminus \iota(D)$. We claim that 
$\lim_{\iota(x) \rightarrow \xi} F(x)$
exists and is in $\partial \Omega$. Fix a sequence $(x_n)_{n\in \Nb}\subset D$ so that $\iota(x_n) \rightarrow \xi$. Since $\overline{\Omega}$ is compact we can assume, passing to a subsequence if necessary, that $F(x_n)$ converges to some $\eta \in \overline{\Omega}$. Fix a point $x_0\in D$. Since $\iota(x_n) \rightarrow \xi \in D^* \setminus \iota(D)$ and $(D, K_D)$ is proper we see that 
\begin{align*}
\lim_{n \rightarrow \infty} K_D(x_n, x_0) = \infty.
\end{align*}
Then, since $F$ is a quasi-isometric embedding, 
\begin{align*}
\lim_{n \rightarrow \infty} K_\Omega(F(x_n), F(x_0)) = \infty.
\end{align*}
Thus $\eta \in \partial \Omega$. Now we claim that 
\begin{align*}
 \lim_{\iota(x) \rightarrow \xi} F(x) = \eta.
\end{align*}
If not, then we would have a sequence $y_n \in D$ so that $\iota(y_n) \rightarrow \xi$, $F(y_n) \rightarrow \eta^\prime$, and $\eta \neq \eta^\prime$. Let $\sigma_n:[0,T_n] \rightarrow D$ be a geodesic with $\sigma_n(0)=x_n$ and $\sigma_n(T_n) = y_n$. Then $(F \circ \sigma_n) : [0,T_n] \rightarrow \Omega$ is a quasi-geodesic and since $\eta \neq \eta^\prime$, Proposition~\ref{prop:qg_visibility} implies that 
\begin{align*}  
\max_{n \in \Nb} K_\Omega(F(x_0), F \circ \sigma_n) < \infty.
\end{align*}
But since $F$ is a quasi-isometric embedding this implies that 
\begin{align*}  
\max_{n \in \Nb} K_D(x_0, \sigma_n) < \infty.
\end{align*}
This contradicts the fact that $(\iota, D^*)$ is a good compactification. Thus for any $\xi \in D^* \setminus \iota(D)$ 
\begin{align*}
 \lim_{\iota(x) \rightarrow \xi} F(x)
\end{align*}
exists and is in $\partial\Omega$.

Next define the map $\wt{F}: D^* \rightarrow \overline{\Omega}$ by 
\begin{align*}
 \wt{F}(\xi) = \begin{cases}
 			F(\iota^{-1}(\xi)), &\text{if $\xi \in \iota(D)$}, \\
 			\lim_{\iota(x) \rightarrow \xi} F(x), &\text{if $\xi \in D^* \setminus \iota(D)$}.
 			\end{cases}
\end{align*}
We claim that $\wt{F}$ is continuous on $D^*$. Since $F$ is continuous on $D$ and $\iota(D) \subset D^*$ is open, it is enough to show that $\wt{F}$ is continuous at each $\xi \in D^* \setminus \iota(D)$. So fix some $\xi \in D^*\setminus \iota(D)$. Since $\overline{\Omega}$ is compact it is enough to show the following:  if $\xi_n \rightarrow \xi$ and $\wt{F}(\xi_n) \rightarrow \eta$ then  $\eta = \wt{F}(\xi) $. Now for each $n$ pick $x_n \in D$ sufficiently close to $\xi_n$ (in the topology of $D^*$) so that $\iota(x_n) \rightarrow \xi$ and 
\begin{align*}
 \|F(x_n) - \wt{F}(\xi_n)\| < 1/n.
\end{align*}
Then 
\begin{align*}
\eta = \lim_{n \rightarrow \infty} \wt{F}(\xi_n) =\lim_{n \rightarrow \infty} F(x_n) 
\end{align*}
but since $\iota(x_n) \rightarrow \xi$, from the discussion in the preceding paragraph, we get 
\begin{align*}
  \lim_{n \rightarrow \infty} F(x_n) = \wt{F}(\xi).
\end{align*}
Hence $\wt{F}$ is continuous. 
\end{proof}

\subsection{The behavior of the Gromov product on Goldilocks domains}\label{ssec:Gromov_product}
Returning to the discussion at the start of this section: if $(X,d)$ is a proper Gromov hyperbolic metric space and $x_n, y_m$ are two sequences in $X$ converging to distinct points in $X(\infty)$ then (by definition)
\begin{align*}
\limsup_{n,m \rightarrow \infty} (x_n | y_m)_{o} < \infty
\end{align*}
for any $o \in X$. We will now show that the Kobayashi distance on a Goldilocks domain has similar behavior. If $\Omega \subset \Cb^d$ is a domain and $x,y,o \in \Omega$, we shall denote the Gromov product for $(\Omega, K_{\Omega})$ by $(x|y)_o^{\Omega}$.

\begin{proposition}\label{prop:gromov_prod}
Suppose $\Omega \subset \Cb^d$ is a Goldilocks domain. If $x_n, y_n \in \Omega$, $x_n \rightarrow \xi \in \partial \Omega$, $y_n \rightarrow \eta \in \partial \Omega$, and $\xi \neq \eta$ then 
\begin{align*}
\limsup_{n,m \rightarrow \infty} (x_n|y_m)_o^{\Omega} < \infty
\end{align*}
for any $o \in \Omega$. 
\end{proposition}

This proposition follows immediately from the next lemma, Proposition~\ref{prop:almost_geod_exist}, and Theorem~\ref{thm:visible}.

\begin{lemma}
Suppose $\Omega \subset \Cb^d$ is a domain and $x,y,o \in \Omega$. If $\sigma:[0,T] \rightarrow \Omega$ is an $(1,\kappa)$-almost-geodesic with $\sigma(0)=x$ and $\sigma(T)=y$ then 
\begin{align*}
(x|y)_o^{\Omega} \leq \frac{3}{2}  \kappa + K_\Omega(o,\sigma).
\end{align*}
\end{lemma}

\begin{proof}
Suppose $s\in[0,T]$ then 
\begin{align}
K_\Omega(x,y) 
&\geq \abs{T-0} - \kappa = \abs{T-s}+\abs{s-0} - \kappa \notag \\
& \geq  K_\Omega(x,\sigma(s)) + K_\Omega(\sigma(s),y) - 3\kappa \label{eq:inv_trngl}
\end{align}
so
\begin{align*}
(x|y)_o^{\Omega} 
&\leq \frac{3}{2} \kappa + \frac{1}{2} \left( K_\Omega(x,o) + K_\Omega(o,y) - K_\Omega(x,\sigma(s))-K_\Omega(\sigma(s),y)\right)\\ 
&\leq \frac{3}{2}  \kappa+ K_\Omega(o, \sigma(s))
\end{align*}
by the reverse triangle inequality.
\end{proof}

\section{Proper holomorphic maps}\label{sec:proper}

The main result of this section once more highlights the point\,---\,but in a manner different from that illustrated by subsection~\ref{ssec:all_about_qi}\,---\,that the conditions defining a Goldilocks domain $\Omega\Subset \Cb^d$ impose adequate control on the oscillation of the values of a proper map into $\Omega$ along suitably chosen sequences approaching the boundary.

Since proper holomorphic maps are, in general, rather far from quasi-isometries of the Kobayashi distance, the methods in this section differ from those in Section~\ref{sec:gromov_prod}. This is also the reason that the statement of Theorem~\ref{thm:proper} addresses a subclass of the class of Goldilocks domains.  

We will need the following results.

\begin{result}[a paraphrasing of Theorem~1 of \cite{DF1977}]\label{res:diederich_fornaess}
Let $\Omega\subset \Cb^d$ be a bounded pseudoconvex domain with $C^2$-smooth boundary. Then there is a defining function $\rho$ of class $C^2$  and a number $\eta_0\in (0, 1)$ such that for each $\eta$, $0< \eta\leq \eta_0$, the function $\wh{\rho} := -(-\rho)^\eta$ is a bounded strictly plurisubharmonic exhaustion function on $\Omega$.
\end{result}

The next result is a version of a Hopf lemma for subharmonic functions. This version is Proposition~1.4 of \cite{M1993b}.

\begin{result}\label{res:Hopf_lemma}
Let $\Omega\subset \Cb^d$ be a bounded domain that satisfies an interior-cone condition with aperture $\theta$. Let $\psi : \Omega\rightarrow [-\infty, 0)$ be a plurisubharmonic function. Then, there exist constants $c > 0$ and $\alpha > 1$ ($\alpha = \pi/\theta$) such that
\begin{align*}
\psi(z)\leq -c(\delta_{\Omega}(z))^\alpha
\end{align*}
for all $z\in \Omega$.
\end{result}

The idea of using the Kobayashi metric to study the boundary behavior of proper holomorphic maps goes back to Diederich and Forn{\ae}ss; see \cite{DF1979}. We adapt their idea to maps for which the target space may have non-smooth boundary.

\begin{proof}[The proof of Theorem~\ref{thm:proper}]
By Result~\ref{res:diederich_fornaess}, we can find a $C^2$-smooth defining function $\rho$ of $D$ and an $\eta\in (0, 1)$ such that $\varphi(z) := -(-\rho(z))^\eta$, $z\in D$, is strictly plurisubharmonic on $D$. Define
\begin{align*}
\psi(w) := \max\left\{\varphi(z) : F(z) = w\right\}
\end{align*}
for each $w\in \Omega$. The function $\psi$, being locally plurisubharmonic at each point not in the branch locus of $F$, is plurisubharmonic away from the branch locus of $F$. As $\psi$ is continuous on $\Omega$, it follows from classical facts\,---\,see, for instance, \cite[Appendix~PSH]{JP1993}\,---\,that $\psi$ is a strictly negative plurisubharmonic function on $\Omega$. As $\Omega$ satisfies an interior-cone condition, there exists, by Result~\ref{res:Hopf_lemma}, a $c > 0$ and an $\alpha > 1$ such that
\begin{align*}
\psi(w)\leq -c(\delta_{\Omega}(w))^\alpha
\end{align*}
for all $w\in \Omega$. Hence
\begin{align}
 (\delta_{\Omega}(F(z)))^\alpha\leq \frac{1}{c}|\psi(F(z))|\leq \frac{1}{c}|\varphi(z)|
 &= \frac{1}{c}|\rho(z)|^\eta \notag \\
 &\leq C\delta_D(z)^\eta \; \; \; \text{for all } z\in D, \label{eq:dist_comparison}
\end{align}
for some $C > 0$, where the last inequality follows from the fact that $\rho$ is a defining function.

It follows from the \emph{proof} of part~(2) of Proposition~\ref{prop:lip}\,---\,see the inequality \eqref{eq:k-metric_optimal}\,---\,that
\begin{align*}
k_{D}(z; v) \leq \frac{\|v\|}{\delta_D(z)}
\end{align*}
for all $z\in D$ and $v\in \Cb^d$. Fix a vector $v$ such that $\|v\| = 1$. Then,
\begin{align*}
k_{\Omega}\left(F(z); F^\prime(z)v\right)
\leq k_D(z; v) \leq \frac{1}{\delta_D(z)}
\end{align*}
for all $z\in D$. It follows from this and \eqref{eq:dist_comparison} that
\begin{align}
 \|F^\prime(z)v\| \leq \frac{(\delta_D(z))^{-1}}{k_{\Omega}\!\left(F(z); 
  \tfrac{F^\prime(z)v}{\|F^\prime(z)v\|}\right)} \leq \
 &\delta_D(z)^{-1}M_{\Omega}(C(\delta_D(z))^{\eta/\alpha}) \notag \\
 &\forall z\in D \text{ and } \forall v\notin {\rm Ker}(F^\prime(z)) : \|v\| = 1, \label{eq:deriv_bd}
\end{align}
and, clearly, the bound on $\|F^\prime(z)v\|$ extends trivially to all unit vectors in ${\rm Ker}(F^\prime(z))$.

As $D$ is bounded and has $C^2$ boundary, there exists an $R > 0$ such that
\begin{align*}
\{z\in D : \delta_D(z)\leq R\}\cup \partial{D} = \bigsqcup_{\xi\in \partial{D}}\{\xi+t\nrml(\xi) : 0\leq t\leq R\},
\end{align*}
where $\nrml(\xi)$ is the inward unit normal vector to $\partial{D}$ at $\xi$. By construction, for each $r\in (0, R]$, we have homeomorphisms $\pi_r : \partial{D}\rightarrow \{z\in D : \delta_D(z) = r\} =: \partial{D}_r$ defined as
\begin{align*}
 \pi_r(\xi) :=&\ \text{the unique $z\in \partial{D}_r$ such that $\delta_D(z) = d_{{\rm Euc}}(\xi, z)$} \\
 =&\ \xi + r\nrml(\xi).
\end{align*}
Pick and fix an $r\in (0, R)$. Write $F = (F_1,\dots, F_d)$ and fix a $j : 1\leq j\leq d$. If $\xi\in \partial{D}$ and $0 < t < r$, then
\begin{align*}
F_j(\xi + t\nrml(\xi)) =
F_j(\pi_r(\xi)) - \int_t^r\left[\sum_{l=1}^{d}\partial_lF_j(\xi + s\nrml(\xi))\nrml(\xi)_l\right]ds,
\end{align*}
where $\partial_l$ denotes the complex differential operator $\partial/\partial z_l$. By \eqref{eq:deriv_bd} and the sentence following it, we have
\begin{align*}
 \int_t^r\left|\sum_{l=1}^{d}\partial_lF_j(\xi + s\nrml(\xi))\nrml(\xi)_l\right|ds
 &\leq \int_t^r\frac{M_{\Omega}(Cs^{\eta/\alpha})}{s}ds \\
 &= \frac{\alpha}{\eta}\int_{Ct^{\eta/\alpha}}^{Cr^{\eta/\alpha}}\frac{M_{\Omega}(u)}{u}du.
\end{align*}
Thus, given that $u\longmapsto M_{\Omega}(n)/u$ is integrable on $[0, R]$, the limit
\begin{align*}
\bv{F}_j(\xi)\,:=\,F_j(\pi_r(\xi)) -
\lim_{t\to 0^+}\int_t^r\left[\sum_{l=1}^{d}\partial_lF_j(\xi + s\nrml(\xi))\nrml(\xi)_l\right]ds
\end{align*}
exists for every $\xi\in \partial{D}$.

We shall now use an aspect of a Hardy--Littlewood-type argument to complete the proof. Pick an $\epsilon > 0$. The preceding argument shows that
\begin{equation}\label{eq:unif_small}
|\bv{F}_j(\xi) - F_j(\pi_r(\xi))| \leq
\frac{\alpha}{\eta}\int_{0}^{Cr^{\eta/\alpha}}\frac{M_{\Omega}(u)}{u}du \; \; \;
\forall \xi\in \partial{D} \text{ and } \forall r\in (0, R).
\end{equation}
Hence, as $u\longmapsto M_{\Omega}(n)/u$ is integrable on $[0, R]$, given $\xi_1, \xi\in \partial{D}$, we can find a constant $r(\epsilon) > 0$ sufficiently small that
\begin{equation}\label{eq:radial_est}
|\bv{F}_j(\xi_i) - F_j(\pi_{r(\epsilon)}(\xi_i))| < \epsilon/3, \; \; \; i = 1, 2.
\end{equation}
Now, as $\left(\left.F_j\right|_{\partial{D}_{r(\epsilon)}}\right)\circ \pi_{r(\epsilon)}$ is uniformly continuous, $\exists\delta > 0$ such that
\begin{align*}
|F_j(\pi_{r(\epsilon)}(\xi_1)) - F_j(\pi_{r(\epsilon)}(\xi_2))| < \epsilon/3\,\text{ whenever }
d_{{\rm Euc}}(\xi_1, \xi_2) < \delta.
\end{align*}
From this and \eqref{eq:radial_est}, we deduce that $\bv{F}_j$ is continuous.

Now write
\begin{align*}
\wt{F}(z) = (\wt{F}_1,\dots, \wt{F}_d)(z)
= \begin{cases}
    F(z), &\text{if $z\in D$}, \\
    \bv{F}(z), &\text{if $z\in \partial{D}$}.
   \end{cases}
\end{align*}
To prove that $\wt{F}$ is continuous on $\conj{D}$, it suffices to show that given a $\xi\in \partial{D}$ and any sequence $\{z_n\}\subset \conj{D}\setminus\{\xi\}$ that converges to $\xi$, $\wt{F}_j(z_n)\to \bv{F}_j(\xi)$ for each $j = 1,\dots, d$. We will construct an auxiliary sequence in $\conj{D}\setminus\{\xi\}$. To this end, consider the continuous map $\pr : (\{z\in D : \delta_D(z)\leq R\}\cup \partial{D})\rightarrow \partial{D}$, defined as
\begin{equation*}
\pr(z) = \pi^{-1}_r(z) \; \; \; \text{if $z\in \partial{D}_r$}.
\end{equation*}
For all {\em sufficiently large} $n$, we can define
\begin{align*}
 Z_n
 := \begin{cases}
       z_n, &\text{if $z_n\in \partial{D}$}, \\
       z_n, &\text{if $z_n\in \{\xi + t\nrml(\xi) : 0 < t\leq R\}$}, \\
       \pr(z_n), &\text{otherwise}.
      \end{cases}
\end{align*}
By continuity of $\pr$, $Z_n\to \xi$. Using integrability of $u\longmapsto M_{\Omega}(n)/u$ once again, it follows from \eqref{eq:unif_small} that $(\wt{F}_j(z_n) - \wt{F}_j(Z_n))\to 0$. However, it follows from the previous two paragraphs that $\wt{F}_j(Z_n)\to \bv{F}_j(\xi)$ for each $j = 1,\dots, d$.  Hence, by the preceding discussion, we infer that $\wt{F}$ is continuous.
\end{proof}

\section{Wolff--Denjoy theorems}\label{sec:WD}

Before proving the Wolff--Denjoy theorems stated in the introduction, let us explain the main idea. Suppose that $\Omega \subset \Cb^d$ is a Goldilocks domain and $f:\Omega \rightarrow \Omega$ is 1-Lipschitz with respect to the Kobayashi metric. The difficult case to rule out is when there exist two sequences $m_i, n_j \rightarrow \infty$ so that $f^{m_i}(o) \rightarrow \xi \in \partial \Omega$, $f^{n_j}(o) \rightarrow \eta \in \partial \Omega$, and $\xi \neq \eta$. In this case we will obtain a contradiction by considering $K_\Omega(f^{m_i}(o), f^{n_j}(o))$. If we assume that $m_i > n_j$ then 
\begin{align*}
 K_\Omega(f^{m_i}(o), f^{n_j}(o)) \leq K_\Omega(f^{m_i - n_j}(o), o).
\end{align*}
Now, if $i \gg j$ then $f^{m_i - n_j}(o)$ should be close to $\xi$. In particular, $ K_\Omega(f^{m_i}(o), f^{n_j}(o))$ is bounded by the ``distance'' from $o$ to $\xi$. On the other hand the visibility condition tells us that any length minimizing curve joining $f^{m_i}(o)$ to $f^{n_j}(o)$ has to pass close to $o$ and so 
\begin{align*}
 K_\Omega(f^{m_i}(o), f^{n_j}(o)) \approx K_\Omega(f^{m_i}(o), o)+K_\Omega(o, f^{n_j}(o))
\end{align*}
which for $i, j \gg 0$ is roughly the sum of the ``distance'' from $o$ to $\xi$  and the ``distance'' from $o$ to $\eta$. Combining these two observations gives a contradiction.

To obtain the second estimate we will use the following observation:

\begin{lemma}\label{lem:middle_pt}
Suppose $\Omega \subset \Cb^d$ is a bounded domain. If $\sigma:[a,b] \rightarrow \Omega$ is a $(1,\kappa)$-quasi-geodesic then for all $t \in [a,b]$ we have 
\begin{align*}
K_\Omega(\sigma(a),\sigma(b)) \leq K_\Omega(\sigma(a),\sigma(t)) + K_\Omega(\sigma(t), \sigma(b)) \leq K_\Omega(\sigma(a),\sigma(b)) + 3\kappa.
\end{align*}
\end{lemma}

\begin{proof}
This follows immediately from the triangle inequality and the definition of a quasi-geodesic. \end{proof}

\subsection{The metric case}
In this section, we give the proof of Theorem~\ref{thm:m_WD}. This theorem is the consequence of Theorem~\ref{thm:metric_WD}, which we now prove. The proof of Theorem~\ref{thm:metric_WD} uses our visibility result and an argument from a paper of Karlsson~\cite[Theorem 3.4]{K2001} about the iterations of 1-Lipschitz maps on general metric spaces.

\begin{theorem}\label{thm:metric_WD}
Suppose $\Omega \subset \Cb^d$ is a Goldilocks domain. If $f:\Omega \rightarrow \Omega$ is 1-Lipschitz with respect to the Kobayashi distance and 
\begin{equation*}
\lim_{n \rightarrow \infty} K_\Omega( f^n(o), o) = \infty
\end{equation*}
for some (hence any) $o \in \Omega$, then there exists a $\xi \in \partial \Omega$ such that 
\begin{align*}
 \lim_{k \rightarrow \infty} f^{k}(x) =\xi
\end{align*}
for all $x \in \Omega$. 
\end{theorem}

\begin{proof}
Fix $o \in \Omega$ and pick a subsequence $m_i  \rightarrow \infty$ so that 
\begin{align*}
K_{\Omega}(f^{m_i}(o), o) \geq K_{\Omega}(f^{n}(o), o)
\end{align*}
for all $n \leq m_i$. By passing to another subsequence we may suppose that $f^{m_i}(o) \rightarrow \xi \in \partial \Omega$. 

Suppose that $f^{n_j}(x) \rightarrow \eta$ for some $x \in \Omega$ and sequence $n_j\rightarrow \infty$. We claim that $\eta=\xi$. Pick a sequence $i_j \rightarrow \infty$ with $m_{i_j} > n_j$.  Now let $\sigma_j :[0,T_j] \rightarrow \Omega$ be an $(1,1)$-almost-geodesic with $\sigma_j(0) = f^{m_{i_j}}(o)$ and $\sigma_j(T_j) = f^{n_j}(x)$. Since $f^{m_i}(o) \rightarrow \xi$, $f^{n_j}(x) \rightarrow \eta$, and $\xi \neq \eta$, Theorem~\ref{thm:visible} implies the existence of some $R > 0$ so that 
\begin{align*}
\max_{j \in \Nb} K_\Omega(o, \sigma_j) \leq R.
\end{align*}
So pick some $t_j \in [0,T_j]$ with
\begin{align*}
K_\Omega(o, \sigma_j(t_j)) \leq R.
\end{align*}
Then by Lemma~\ref{lem:middle_pt} we have 
\begin{align*}
K_{\Omega}(f^{m_{i_j}}(o), f^{n_j}(x))
& \geq K_{\Omega}(f^{m_{i_j}}(o),\sigma_j(t_j)) + K_\Omega(\sigma_j(t_j), f^{n_j}(x)) -3 \\
& \geq K_\Omega(f^{m_{i_j}}(o), o) + K_\Omega(o, f^{n_j}(x))- 3-2R
\end{align*}
On the other hand
\begin{align*}
K_{\Omega}(f^{m_{i_j}}(o), f^{n_j}(x)) \leq K_{\Omega}(f^{m_{i_j}-n_j}(o), o)+K_\Omega(o,p) \leq K_\Omega(f^{m_{i_j}}(o), o)+K_\Omega(o,x).
\end{align*}
So 
\begin{align*}
K_\Omega(o, f^{n_j}(x)) \leq 3+2R+K_\Omega(o,x)
\end{align*}
and we have a contradiction.
\end{proof}

Finally, we provide

\begin{proof}[The proof of Theorem~\ref{thm:m_WD}]
 Since $(\Omega, K_\Omega)$ is Cauchy complete, a result of Ca{\l}ka~\cite[Theorem 5.6]{C1984b} implies that 
either 
 \begin{align*}
 \lim_{n \rightarrow \infty} K_{\Omega}(f^n(x), x) = \infty
 \end{align*}
for any $x \in \Omega$ or
 \begin{align*}
 \sup_{n \geq 0} K_{\Omega}(f^n(x), x) < \infty
 \end{align*}
 for any $x \in \Omega$. In the first case, Theorem~\ref{thm:metric_WD} implies that there exists $\xi \in \partial \Omega$ so that
\begin{equation*}
\lim_{n \rightarrow \infty} f^n(x) = \xi
\end{equation*}
 for any $x \in \Omega$. In the second case, Result~\ref{res:hopf_rinow} implies that the orbit $\{ f^n(x): n \in \Nb\}$ is relatively compact in $\Omega$ for any $x \in \Omega$.
\end{proof}

\subsection{The holomorphic case}

We shall now give a proof of Theorem~\ref{thm:WD}.

\begin{lemma}\label{lem:limits}
Let $\Omega \subset \Cb^d$ be a Goldilocks domain. Suppose $f:\Omega \rightarrow \Omega$ is a holomorphic map. If $f^{n_i}$ converges to some $F: \Omega \rightarrow \partial \Omega$ then $F \equiv \xi$ for some $\xi \in \partial \Omega$.
\end{lemma}

\begin{proof} Fix some $x \in \Omega$. Then $\lim_{i \rightarrow \infty} d(f^{n_i})_x=dF_x$. And if $v \in \Cb^d$ then 
\begin{align*}
k_\Omega(f^{n_i}(x); d(f^{n_i})_x v) \leq k_\Omega(x; v).
\end{align*}
Let $\tau = \max \{ k_\Omega(x; v) : \norm{v} =1\}$. We claim that 
\begin{align*}
\norm{ d(f^{n_i})_x v}   \leq  \tau M_\Omega(\delta_\Omega(f^{n_i}(x))) \text{ when } \norm{v}=1. 
\end{align*}
It clearly suffices to consider the case when $d(f^{n_i})_x v \neq 0$. In this case
\begin{align*}
1 \leq \frac{k_\Omega(x; v)}{k_\Omega(f^{n_i}(x); d(f^{n_i})_x v)} \leq \frac{\tau}{k_\Omega(f^{n_i}(x); d(f^{n_i})_x v)}.
\end{align*}
Then 
\begin{align*}
\norm{ d(f^{n_i})_x v}  
& \leq \frac{\tau\norm{ d(f^{n_i})_x v}}{k_\Omega(f^{n_i}(x); d(f^{n_i})_x v)} = \frac{\tau}{k_\Omega\left(f^{n_i}(x); \frac{d(f^{n_i})_x v}{\norm{d(f^{n_i})_x v}}\right)} \\
& \leq \tau M_\Omega(\delta_\Omega(f^{n_i}(x))). 
\end{align*} 

Then since $\delta_\Omega(f^{n_i}(x)) \rightarrow 0$ and $\lim_{i \rightarrow \infty} d(f^{n_i})_x=dF_x$ we see that $dF_x=0$. Since $x \in \Omega$ was arbitrary we see that $dF=0$ and hence $F$ is constant.
\end{proof}

\begin{proof}[The proof of Theorem~\ref{thm:WD}]
Since $\Omega$ is taut by~\cite[Theorem 2.4.3]{A1989}, either
\begin{enumerate}
\item for any $x \in \Omega$, the orbit $\{ f^n(x): n \in \Nb\}$ is relatively compact in $\Omega$; or
\item for any $x \in \Omega$, 
\begin{equation*}
\lim_{n \rightarrow \infty} d_{\Euc}(f^n(x), \partial \Omega) = 0.
\end{equation*}
\end{enumerate}
Suppose that the second condition holds. Montel's theorem tells us that there exist subsequences $\{f^{n_j}\}$ that converge locally uniformly to $\partial\Omega$-valued functions. By Lemma~\ref{lem:limits}, the latter are constant functions. Thus, we will identify the set
\begin{align*}
 \Gamma:=\overline{ \{f^n : n \in \Nb\}}^{{\rm compact-open}}\setminus \{f^n : n \in \Nb\}
\end{align*}
as a set of points in $\partial\Omega$. Our goal is to show that $\Gamma$ is a single point.

Assume for a contradiction that $\Gamma$ is not a single point. 
\medskip

\noindent \textbf{Case 1:} Suppose that for some (hence any) $o \in \Omega$ we have 
\begin{equation*}
\limsup_{n \rightarrow \infty} K_\Omega(f^n(o), o) = \infty.
\end{equation*}
Then we can find a subsequence $m_i \rightarrow \infty$ so that 
\begin{align*}
K_\Omega(f^{m_i}(o), o) \geq &\;K_\Omega(f^k(o), o) \; \; \text{for all } k\leq m_i.
\end{align*}
By passing to a subsequence we can assume that $f^{m_i} \rightarrow \xi \in \partial \Omega$. Now by assumption, there exists a subsequence $n_j \rightarrow \infty$ so that $f^{n_j} \rightarrow \eta \in \partial \Omega$ and $\eta \neq \xi$.
\medskip

\noindent \textbf{Case 1(a):} First consider the case in which
\begin{equation*}
\limsup_{j \rightarrow \infty} K_\Omega(f^{n_j}(o), o) = \infty.
\end{equation*}
In this case we can repeat the proof of Theorem~\ref{thm:metric_WD} essentially verbatim: Pick $i_j \rightarrow \infty$ so that $m_{i_j} > n_j$. Now let $\sigma_j :[0,T_j] \rightarrow \Omega$ be an $(1,1)$-almost-geodesic with $\sigma_j(0) = f^{m_{i_j}}(o)$ and $\sigma_j(T_j) = f^{n_j}(o)$. Since $f^{m_i} \rightarrow \xi$, $f^{n_j} \rightarrow \eta$, and $\xi \neq \eta$, Theorem~\ref{thm:visible} implies the existence of some $R > 0$ so that 
\begin{align*}
\max_{j \in \Nb} K_\Omega(o, \sigma_j) \leq R.
\end{align*}
So pick some $t_j \in [0,T_j]$ with
\begin{align*}
K_\Omega(o, \sigma_j(t_j)) \leq R.
\end{align*}
Then by Lemma~\ref{lem:middle_pt} we have 
\begin{align*}
K_{\Omega}(f^{m_{i_j}}(o), f^{n_j}(o))
& \geq K_{\Omega}(f^{m_{i_j}}(o),\sigma_j(t_j)) + K_\Omega(\sigma_j(t_j), f^{n_j}(o)) -3 \\
& \geq K_\Omega(f^{m_{i_j}}(o), o) + K_\Omega(o, f^{n_j}(o))- 3-2R
\end{align*}
On the other hand
\begin{align*}
K_{\Omega}(f^{m_{i_j}}(o), f^{n_j}(o)) \leq K_{\Omega}(f^{m_{i_j}-n_j}(o), o) \leq K_\Omega(f^{m_{i_j}}(o), o).
\end{align*}
So 
\begin{align*}
K_\Omega(o, f^{n_j}(o)) \leq 3+2R
\end{align*}
and we have a contradiction.
\medskip

\noindent \textbf{Case 1(b):} Next consider the case in which
\begin{equation*}
\limsup_{j \rightarrow \infty} K_\Omega(f^{n_j}(o), o) < \infty.
\end{equation*}
By Lemma~\ref{lem:limits}, for any $\ell \in \Nb$ we have
\begin{align*}
 \lim_{j \rightarrow \infty} f^{n_j-\ell}(o) = \eta.
\end{align*}
Let 
\begin{align*}
M_{\ell}:= \limsup_{j \rightarrow \infty} K_\Omega(f^{n_j-\ell}(o), o).
\end{align*}
We claim that 
\begin{align*}
\limsup_{\ell \rightarrow \infty} M_{\ell} < \infty.
\end{align*}
Suppose not; then we find $\ell_k \rightarrow \infty$ so that $M_{\ell_k} > k$,
$k = 1, 2, 3,\dots$ Then we pick $j_k \rightarrow \infty$ so that 
\begin{align*}
K_\Omega(f^{n_{j_k}-\ell_k}(o), o) > k, \; \; \text{and} \; \; d_{\Euc}(f^{n_{j_k}-\ell_k}(o), \eta) < 1/k.
\end{align*}
But then $f^{n_{j_k}-\ell_k}(o) \rightarrow \eta$ and 
\begin{align*}
\lim_{k \rightarrow \infty} K_\Omega(f^{n_{j_k}-\ell_k}(o), o) =\infty
\end{align*}
which is impossible by Case 1(a). So we see that 
\begin{align*}
\limsup_{\ell \rightarrow \infty} M_{\ell} < \infty.
\end{align*}

Then
\begin{align*}
\limsup_{i \rightarrow \infty} \limsup_{j \rightarrow \infty} K_\Omega(f^{m_i} (o), f^{n_j} (o)) &\leq \limsup_{i \rightarrow \infty} \limsup_{j \rightarrow \infty} K_\Omega( o, f^{n_j-m_i} (o))\\
&= \limsup_{i \rightarrow \infty} M_{m_i} < \infty,
\end{align*}
and 
\begin{align*}
\limsup_{i \rightarrow \infty} \limsup_{j \rightarrow \infty} & K_\Omega(f^{m_i} (o), f^{n_j}( o))\\
& \geq \limsup_{i \rightarrow \infty} \limsup_{j \rightarrow \infty}\Big(K_\Omega(f^{m_i}(o),  o) - K_\Omega(o, f^{n_j}(o))\Big) \\
& \geq \limsup_{i \rightarrow \infty} \Big(K_\Omega(f^{m_i} (o),o) - M_{0} \Big)=\infty. 
 \end{align*}
 So we again have a contradiction.
 \medskip
 
 \noindent \textbf{Case 2:} Suppose that for some (hence any) $o \in \Omega$ we have 
\begin{equation*}
\limsup_{n \rightarrow \infty} K_\Omega(f^n(o), o) <\infty.
\end{equation*}

Suppose that $\xi, \eta \in \Gamma$ are two distinct points. Fix neighborhoods $V_\xi$ of $\xi$ and $V_\eta$ of $\eta$
so that $\overline{V_\xi} \cap \overline{V_\eta} = \emptyset$. By Theorem~\ref{thm:visible} there exists a
compact set $K \subset \Omega$ with the following property: if $\sigma:[0,T] \rightarrow \Omega$ is any $(1,2)$-almost-geodesic satisfying $\sigma(0) \in V_\xi$ and $\sigma(T) \in V_\eta$ then $\sigma([0,T]) \cap K \neq \emptyset$.
\vspace{1mm}

Next, for $\delta > 0$ define the function $G_\delta: K \times K \rightarrow \Rb$ by 
\begin{align*}
G_\delta(k_1, k_2) := \inf\{ K_\Omega(f^m (k_1), k_2) : d_{\Euc}(f^m(k_1),\xi) < \delta\}.
\end{align*}
By the assumptions for Case~2,
\begin{align*}
\sup\{ G_\delta(k_1, k_2) : \delta >0 , k_1, k_2 \in K\} < \infty
\end{align*}
and if $\delta_1 < \delta_2$ then $G_{\delta_1} \geq G_{\delta_2}$. So the function
\begin{align*}
G(k_1, k_2) := \lim_{\delta \rightarrow 0} G_\delta(k_1, k_2)
\end{align*}
is well defined.

Next, let
 \begin{align*}
 \epsilon: = \liminf_{z \rightarrow \eta} \inf_{k \in K} K_\Omega(k, z).
 \end{align*}
By Proposition~\ref{prop:lip},  $\epsilon > 0$. Now pick $q_1, q_2 \in K$ so that 
\begin{align*}
G(q_1, q_2) < \epsilon + \inf\{ G(k_1, k_2) : k_1, k_2 \in K\}.
\end{align*}

Fix a sequence of integers $n_j \rightarrow \infty$ so that $f^{n_j} \rightarrow \eta$. Suppose $\mu_i \rightarrow \infty$ is any sequence of integers such that $f^{\mu_i} \rightarrow \xi$. Then by Lemma~\ref{lem:limits}
\begin{align*}
 \lim_{i \rightarrow \infty} f^{\mu_i+n_j}(o) = \lim_{i \rightarrow \infty} f^{\mu_i}(f^{n_j}(o)) = \xi.
\end{align*}
So we can find a subsequence $\{\mu_{i_j}\}\subset \{\mu_i\}$ so that $f^{\mu_{i_j} + n_j} \rightarrow \xi$. Therefore, we can find a sequence of integers $m_j \rightarrow \infty$ such that
\begin{align*}
f^{m_j} &\rightarrow \xi, \\
f^{m_j+n_j} &\rightarrow \xi, \\
\lim_{j \rightarrow \infty} K_\Omega( f^{m_j}(q_1), q_2) &= G(q_1, q_2).
\end{align*}
Finally, fix a sequence $\kappa_j \searrow 0$ with $\kappa_j \leq 2$. By Proposition~\ref{prop:almost_geod_exist}, there exists an $(1,\kappa_j)$-almost-geodesic $\sigma_{j}: [0,T_j] \rightarrow \Omega$ with $\sigma(0)=f^{m_j + n_j}(q_1)$ and $\sigma(T_j) =f^{n_j}(q_2)$. For $j$ sufficiently large, $\sigma_j(0) \in V_\xi$ and $\sigma_j(T_j) \in V_\eta$. Since each $\sigma_j$ is an $(1,2)$-almost-geodesic, by the construction of $K$ there exists, for each $j$ sufficiently large, a point $k_j \in K \cap \sigma([0,T_j])$. Then, by Lemma~\ref{lem:middle_pt}, we have
\begin{align*}
K_\Omega(f^{m_j + n_j}(q_1), f^{n_j}(q_2))
\geq K_\Omega(f^{m_j + n_j}(q_1), k_j) + K_\Omega(k_j, f^{n_j}(q_2)) - 3 \kappa_j.
\end{align*}
Now by our definition of $\epsilon$ we have
\begin{align*}
\liminf_{j \rightarrow \infty} K_\Omega(k_j, f^{n_j}(q_2)) \geq \epsilon. 
\end{align*}
After passing to a subsequence we can suppose that $k_j \rightarrow k \in K$. Then since  $f^{m_j + n_j}(q_1) \rightarrow \xi$ we see that 
\begin{align*}
\liminf_{j \rightarrow \infty} & K_\Omega(f^{m_j + n_j}(q_1), k_j) \geq
\liminf_{j \rightarrow \infty} \Big( K_\Omega(f^{m_j + n_j}(q_1), k)- K_\Omega(k,k_j) \Big) \\
& = \liminf_{j \rightarrow \infty} K_\Omega(f^{m_j + n_j}(q_1), k) \geq G(q_1, k).
\end{align*}
Since $ \kappa_j \rightarrow 0$, from the last three estimates, we get 
\begin{align*}
\liminf_{j \rightarrow \infty}  K_\Omega(f^{m_j + n_j}(q_1), f^{n_j}(q_2)) \geq G(q_1, k) + \epsilon.
\end{align*}
On the other hand,
\begin{align*}
\limsup_{j \rightarrow \infty} K_\Omega(f^{m_j + n_j}(q_1), f^{n_j}(q_2))
\leq  \limsup_{j \rightarrow \infty}K_\Omega(f^{m_j}(q_1), q_2) = G(q_1, q_2).
\end{align*}
So we have
\begin{align*}
G(q_1,q_2) \geq G(q_1,k)+\epsilon
\end{align*}
which contradicts our choice of $q_1, q_2 \in K$.

In both Cases~1 and 2, we obtain contradictions. Hence, $\Gamma$ contains a single point.
\end{proof}

\bibliographystyle{alpha}
\bibliography{complex_kob}

\end{document}